\documentclass[lettersize,journal]{IEEEtran}

\usepackage[T1]{fontenc}
\usepackage{mathrsfs}
\usepackage{bm}
\usepackage{mathrsfs}
\usepackage{url}
\usepackage{tikz}
\usepackage[colorlinks=true,citecolor=blue,urlcolor=blue]{hyperref}
\usepackage[noend]{algpseudocode}
\usepackage{amsmath,amssymb,amsfonts}
\usepackage{algorithmicx,algorithm}
\usepackage{graphicx}
\usepackage{textcomp}
\usepackage{array}
\usepackage{cases}
\usepackage{booktabs}
\usepackage{subfigure}
\usepackage{mfirstuc}
\usepackage{float}
\usepackage{verbatim}
\usepackage{hyperref}
\usepackage{wrapfig}
\usepackage{makecell}
\usepackage{multirow}
\usepackage{amssymb}

\def\BibTeX{{\rm B\kern-.05em{\sc i\kern-.025em b}\kern-.08em
		T\kern-.1667em\lower.7ex\hbox{E}\kern-.125emX}}

\newtheorem{remark}{Remark}

\newtheorem{corollary}{Corollary}
\newtheorem{theorem}{Theorem}

\begin{document}

\title{An Operator Splitting Scheme for Distributed Optimal Load-side Frequency Control with Nonsmooth Cost Functions\\
}

\author{
	Yifan Wang, \IEEEmembership{Student Member, IEEE}, Shuai Liu, \IEEEmembership{Member, IEEE}, Xianghui Cao, \IEEEmembership{Senior Member, IEEE} and Mo-Yuen Chow, \IEEEmembership{Fellow, IEEE}
	\thanks{Y. Wang and X. Cao are with the School of Automation, Southeast University, Nanjing 210096, China (E-mail: evan@seu.edu.cn; xhcao@seu.edu.cn).}
    \thanks{S. Liu is with the School of Control Science and Engineering,
Shandong University, Jinan 250061, China (e-mail: lius0025@ntu.edu.sg).}
    \thanks{M. Y. Chow is with the UM–SJTU Joint Institute Shanghai Jiao Tong
University, Shanghai 200240, China (e-mail: moyuen.chow@sjtu.edu.cn).}
}

\maketitle

\begin{abstract}
The increasing penetration of renewable energy resources and utilization of energy storage systems pose new challenges in maintaining power system's stability. Specifically, the cost function of regulation no longer remains smooth, which complicates the task of ensuring nominal frequency and power balance, particularly in a distributed manner. This paper proposes a distributed proximal primal-dual (DPPD) algorithm, based on a modified primal-dual dynamics equipped with operator splitting technique, to address this nonsmooth frequency regulation problem on load side. By Lyapunov stability and invariance theory, we prove that DPPD algorithm achieves global asymptotic convergence to the optimal solution that minimizes the nonsmooth regulation cost, while restoring the frequency under constraints such as capacity limits of load and power flow. Finally, we demonstrate the effectiveness and robustness of DPPD algorithm by simulations on the IEEE 39-bus system.
\end{abstract}



\section{Introduction}\label{sec:introduction}
In power systems, frequency control regulates the synchronized frequency around nominal value (50 or 60 Hz) to maintain the stability. Traditionally, a three-layer frequency control framework on generation side is used to implement this task, where droop control (primary control), automatic generation control (AGC, secondary control) and economic dispatch (ED, tertiary control) work in concert on different timescales \cite{3ji,shu}. However, the rapid penetration of renewable energy resources and fluctuating demand pose new challenges to maintain power balance and stability in frequency regulation problem. Consequently, these problems may undermine the safety and economy of power system operation. On the other hand, due to the integration of renewables and the withdrawal of synchronous machines, converters substitute synchronous machines with large inertia on the generation side, leading to a relatively slower response compared with that in the load-side control. Therefore, it is a promising complement to conduct researches on load-side frequency control in the new era of power system.


Currently, the most prevalent methods for load control, particularly in a distributed manner, can be categorized into two groups. The first one, known as the partial primal-dual gradient algorithm, is derived from primal-dual framework by addressing the optimal load control (OLC) problem. It is worth noting that most works within this category focus on a joint problem involving ED on a slow timescale and frequency regulation on a fast timescale. In \cite{zch1}, a feedback controller is designed as the distributed algorithm for primary frequency control, but it is incapable to restore the nominal frequency. To mitigate this weakness, schemes that integrate primary and secondary frequency control are proposed \cite{vfolc,pd_caid,chenxin,hill}. In \cite{vfolc}, virtual phase angle is firstly introduced to constrain the real flows, making it possible to interpret the power system dynamics as part of primal-dual algorithm. Following with this work, Cai et al. decompose the joint problem into two subproblems in slow and fast timescales, and solve them without loss of optimality \cite{pd_caid}. \cite{chenxin} proposes an automatic load control (ALC) algorithm with a global exponential convergence if the capacity limits are ignored. It is worth mentioning that the proposed ALC averts the measurement of power imbalance or injection, that is usually unknown in advance. In \cite{hill}, frequency regulation is executed on both generation and load sides under nonlinear swing dynamics.
The second one is referred to as distributed average proportional integral (DAPI) control which enjoys the ease of execution, yet it is not easy to enforce line and power flow constraints \cite{dapi1,dapi2}. Specifically, Simpson-Porco et al. introduce DAPI controller for frequency regulation and voltage control \cite{dapi1}. Further, the stability analysis is provided under weaker assumptions on objective functions and communication topology in \cite{dapi2}. Obviously, both of those two methods have their advantages and drawbacks.
As complements and improvements, new approaches have emerged, e.g., involving reinforcement learning to enhance the frequency stability maybe caused by inexact measurements \cite{cui}, counting on model predictive control together with fuzzy control to reduce the frequency deviation\cite{mpc_pd,fuzzy}, and coping with higher system dynamics in frequency regulation problem \cite{nima}.

However, there still exists a concern not extensively addressed in the field of OLC problem, i.e., the case of nonsmooth objective function. Traditionally, many works adopt the scheme of projecting the partial gradient dynamics directly on a feasible set for constrained optimization problem \cite{zch1,chow,yefeng,tii}. However, in cases with nonsmoothness, the existence of solution is nontrivial due to the nonconvexity resulting from the projection of a discontinuous dynamics. Furthermore, discontinuous characteristics may induce a violent vibration in dynamic systems that usually require Lipschitz continuity. To tackle this issue, related works in distributed optimization can be classified into three types. For the first type, projection on the tangent cone is employed, but it increases additional computation burden \cite{tangent2}. For the second type, distributed projected output feedback (DPOF) algorithm and distributed derivative feedback (DDF) algorithm are devised in virtue of feedback-based principle \cite{pof,Chen2021DistributedMD}. Subsequently, the DPOF scheme is applied to address nonsmooth OLC problem in \cite{zhaojian2}. However, subgradient is directly used in the second one, and it is not dominant in terms of convergence rate. The third type of approach is based on operator splitting technique, usually referred to as proximal-based algorithm \cite{davis2017three,yuewei,9831961,randprox,Distributed-Optimization-With-Coupling-Constraints-2023TAC}.

Notwithstanding, it is still challenging to design continuous-time smooth algorithm for nonsmooth OLC problem with rigorous convergence analysis, which is exactly the motivation of our work. The differences between our work and \cite{davis2017three,yuewei,9831961,randprox,Distributed-Optimization-With-Coupling-Constraints-2023TAC} lie in two folds: 1) In \cite{davis2017three,randprox,Distributed-Optimization-With-Coupling-Constraints-2023TAC}, nonsmooth optimization problem with set constraint was studied with proximal splitting algorithm in discrete time. In \cite{yuewei}, proximal splitting algorithm in continuous time was proposed for distributed nonsmooth consensus problem. While the nonsmooth OLC problem in this work is subject to heterogenous local set constraints, equality and inequality constraints, which can be viewed as an extension in terms of problem establishment. Moreover, the complex problem further makes the analysis more intertwined; 2) \cite{davis2017three,yuewei,9831961,randprox,Distributed-Optimization-With-Coupling-Constraints-2023TAC} merely develop algorithms from the perspective of optimization. However, nonsmooth OLC problem is hierarchical in three control layers on different timescales. Hence, analyzing the proposed nonsmooth optimization algorithm in a unified framework that integrates optimization algorithm with power system dynamics is more challenging.

In this paper, we develop a distributed proximal primal-dual (DPPD) algorithm for OLC problem, where the cost of load adjustment is extended to be a nonsmooth function. Building upon prior researches \cite{zch1,vfolc,chenxin,zhaojian2}, the design of DPPD algorithm starts from the interpretation of power system dynamics as part of primal-dual framework. Therein, we split the nonsmooth part and utilize the operator splitting technique to circumvent the trouble resulting from the nonsmoothness. To summarize, the main contributions of this paper are listed in the following aspects:

1) {It is the first to develop distributed primal-dual algorithm with proximal approach for nonsmooth OLC problem. Compared with \cite{zhaojian2} that focuses on a similar topic but is based on the DPOF scheme, the proposed DPPD algorithm is derived from a modified primal-dual dynamics equipped with operator splitting technique that leverages the separable structure of objective functions. As a byproduct, it reveals that DPPD is a proportional-derivative (PD) or proportional-integral (PI) variant of those protocols in \cite{zch1,vfolc,chenxin,zhaojian2} following from classical primal-dual dynamics. Therefore, the proposed DPPD algorithm brings an essentially different scheme and convergence analysis for extension to nonsmooth OLC problem and its possibly added benefits offered by the derivative or integral feedback.}


2) {The optimality and global asymptotic convergence of DPPD algorithm are proved via Lyapunov stability and invariance theory without Lipschitz continuity assumption. Under the framework of our analysis, the convergence results of DPPD algorithm can be easily extended to the OLC problem with constrained frequency region and constrainted virtual phase angle.}

3) {In the sense of dynamics of optimization variables, the convergence rate is further analyzed in the case whether capacity limits (i.e., inequality constraints) are considered or not, which is shown as $O(1/\sqrt{t})$.}

The organization of this paper is as follows. In Section \uppercase\expandafter{\romannumeral2}, some basic preliminaries are introduced. Section \uppercase\expandafter{\romannumeral3} formulates the power network model and nonsmooth OLC problem. Section \uppercase\expandafter{\romannumeral4} presents the development of DPPD algorithm. In Section \uppercase\expandafter{\romannumeral5}, the proposed DPPD algorithm is analyzed for its optimality and convergence. The effectiveness of DPPD algorithm is validated via some case studies in section \uppercase\expandafter{\romannumeral6}. Finally, it comes to the conclusion in Section \uppercase\expandafter{\romannumeral7}.
\par \emph{Notations}: Throughout this paper, column vectors with compatible dimensions are denoted by boldface letters. $
\mathrm{col}\{x_1,...,x_m\}$ denotes the stacked column vector, where $x_i\in\mathbb{R}$, $\forall i\in \{1,...,m\}$. $\mathbb{R}$, $\mathbb{R}_{+}$ and $\mathbb{R}^{n\times l}$ are sets of real numbers, non-negative real numbers, and real matrices of dimension $n\times l$, respectively. Identity matrix is denoted by $I$ with proper dimension. We use $\|\cdot\|_1$ and $\|\cdot\|$ for $l_1$ norm and standard Euclidean norm, respectively. $\left\langle\cdot,\cdot\right\rangle$ represents the inner product in the Euclidean space. Let superscript $T$ be the transpose operator and superscript $-1$ be the inverse operator.
\section{Preliminaries}
\subsection{Convex Analysis}
A set $\mathcal{X}\subset \mathbb{R}^n$ is convex if it has $\alpha x_1+(1-\alpha)x_2\in \mathcal{X}$ for any $x_1, x_2\in \mathcal{X}$ and any $\alpha\in[0,1]$. A function $f: \mathbb{R}^n\rightarrow \mathbb{R}$ is convex (resp. strictly convex) if it has $\alpha f(x_1)+(1-\alpha)f(x_2)\geq$ (resp. $>$) $f(\alpha x_1+(1-\alpha)x_2)$. If there exists a positive constant $\beta$ such that
$$f(x_1)\geq f(x_2)+\left\langle\nabla f(x_2),x_1-x_2\right\rangle+\frac{\beta}{2}\|x_1-x_2\|^2,$$ then $f$ is said $\beta$-strongly convex where $\nabla f(x)$ is the gradient of $f$ at $x$. The subdifferential of a convex function $f$ at $x$ is defined as
$\partial f(x)=\{p\in\mathbb{R}^n| f(x_1)-f(x)\geq p^T(x_1-x),\forall x_1\in\mathbb{R}^n \}$.

The tangent cone of a closed convex set $\mathcal{X}$ at $x\in \mathcal{X}$ is defined as
$$  \mathcal{T}_\mathcal{X}(x)\triangleq\{\lim\limits_{k\rightarrow\infty}\frac{x_k-x}{\tau_k}|x_k\in \mathcal{X},\tau_k>0,x_k\rightarrow x,\tau_k\rightarrow0\}.$$
Besides, the normal cone is defined as $$\mathcal{C}_\mathcal{X}(x)\triangleq\{z|\langle z,y-x\rangle\leq 0,\,\forall y\in \mathcal{X}\}.$$

The projection operator over the closed set $\mathcal{X}$ is defined as
$$\mathcal{P}_\mathcal{X}(y)=\arg\min_{x\in \mathcal{X}}\|x-y\|.$$ Note that there exists the following relationship (\hspace{-0.1pt}\cite{nonlinear} Lemma 2.38) between the normal cone and the projection on the set $\mathcal{X}$, i.e.,
\begin{equation}\label{conedefinition}
\mathcal{C}_\mathcal{X}(x)=\{z|\mathcal{P}_\mathcal{X}(x+z)=x\}.
\end{equation}

\subsection{Proximal Operator}
Let $f(x): \mathbb{R}^n\rightarrow \mathbb{R}$ be a lower semicontinuous convex functions, which is efficiently computable. Then, the proximal operator of $f$ at $y\in \mathbb{R}^n$ with parameter $\tau >0$ is
\begin{equation}\label{prox}
  \mathbf{prox}_{\tau f}(y)=\arg\min\limits_{x}\left\{f(x)+\frac{1}{2\tau}\|x-y\|^2\right\}.
\end{equation}
The projection operator is a special case of the proximal operator, which can be verified by letting $f(x)$ in \eqref{prox} be an indicator function of a closed convex set $\mathcal{X}$ such that $f(x)=0$ if $x\in \mathcal{X}$ and $f(x)=+\infty$ otherwise. Then, some useful inequalities regarding projection operators are valid for proximal operators. Besides, it implies from \eqref{prox} that
\begin{equation}\label{3}
   x=\mathbf{prox}_{\tau f}(y)=(I+\tau \partial f)^{-1}(y).
\end{equation}

\section{System Model and Problem Statement}
\subsection{Power network model}
Consider the power network\footnote{{The communication network has the same graph as the power transmission network, but the information flow is undirected, which is a common configuration as in \cite{zch1,vfolc,chenxin,zhaojian2}. The Laplacian matrix of the communication network can be denoted by $CC^T\in\mathbb{R}^{n\times n}$.}} as a directed and connected graph $\mathcal{G}(\mathcal{N},\mathcal{E})$, where $\mathcal{N}=\{1,..
.,n\}$ denotes the bus set and $\mathcal{E}=\{1,...,l\}\subseteq\mathcal{N}\times\mathcal{N}$ denotes the transmission line set. The incidence matrix associated with the graph $\mathcal{G}$ is denoted by $C\in\mathbb{R}^{n\times l}$, where the $(i,e)$-th element $C_{ie}=1$ if bus $i$ is the source in the directed transmission line, i.e., $e=ij\in \mathcal{E}$, $C_{ie}=-1$ if bus $i$ is the end of transmission line, i.e., $e=(j,i)\in \mathcal{E}$, and $C_{ie}=0$ otherwise. Clearly, it has $\mathbf{1}^T C=\mathbf{0}^T$. $CC^T$ is a symmetric Laplacian matrix, and $CBC^T$ is a $B_{ij}$-weighted Laplacian matrix where $B\in\mathbb{R}^{l\times l}$ is a susceptance matrix. In such a power network, we assume that the there are two types of buses, i.e., the generator and load buses which are denoted as sets $\mathcal{N_G}$ and $\mathcal{N_L}$, respectively. Without loss of generality, the buses in $\mathcal{N_G}$ and $\mathcal{N_L}$ are indexed by $1,...,n_g$ and $n_g+1,...n_g+n_l$, respectively. It is assumed that the generator buses are connected with generators and possibly attached loads, while the load buses are only attached with loads.

For ease of formulation, we make the following commonly adopted assumptions \cite{vfolc,nima}:

1) The voltage magnitude at each bus is fixed.

2) The reactive power does not affect the phase angle and the frequency.

3) The transmission lines are lossless and characterized by the susceptances $B_{ij}$ between any two buses $i$ and $j$.

Consequently, the DC power flow model can be modelled as $P_{ij}=B_{ij}(\theta_i-\theta_j)$, where $\theta_{i}$ is the power phase angle on bus $i$. Then, the overall dynamics of the power system is described by a set of linearized swing equations given as follows (referring to \cite{shu} for more details):
\begin{subequations}\label{dynamics}
  \begin{alignat}{4}
    \dot{\theta}_{i}=&\,\omega_i, \forall i\in\mathcal{N} \label{3a}\\
     \dot{P}_{ij}=&\,B_{ij}(\omega_i-\omega_j), \forall ij\in \mathcal{E} \label{3b}\\
     M_i\dot{\omega}_i=&\,P_i^{in}-D_i\omega_i-d_i+\sum_{k:ki\in\mathcal{E}}\hspace{-8pt}P_{ki}-\sum_{j:ij\in\mathcal{E}}\hspace{-8pt}P_{ij}, \forall i\in\mathcal{N_G} \label{3c} \\
    0=&\,P_i^{in}-D_i\omega_i-d_i+\sum_{k:ki\in\mathcal{E}}\hspace{-8pt}P_{ki}-\sum_{j:ij\in\mathcal{E}}\hspace{-8pt}P_{ij}, \forall i\in \mathcal{N_L}  \label{3d}
  \end{alignat}
\end{subequations}
where $\omega_i$ is the frequency deviation from the nominal value, namely 50 Hz or 60 Hz. $M_i$ and $D_i$ are the generator inertia and damping coefficient, respectively. $d_i$ denotes the power demand of controllable load at node $i$. $P_i^{in}$ denotes the local power injection at node $i$ (i.e., local power generation\footnote{Power generation is also uncontrollable on load side. Hence, compared to $d_i$, $P_i^{in}$ is not a control or optimization variable in OLC problem. In this work, it is assumed that $P_i^{in}$ is known in advance, e.g., by prediction.}  minus local power demand of uncontrollable load), which may be positive or negative. To some extent, $P_i^{in}$ can be viewed as a disturbance on the bus $i$. For brevity, the angle phase deviation and the frequency deviation between buses $i$ and $j$ are denoted as $\theta_{ij}$ and $\omega_{ij}$, respectively. The equations in \eqref{dynamics} present a standard model for the frequency regulation problem, which approximately captures the dynamic characteristics of phase angle, active power flow in branches and frequency \cite{shu}.
\subsection{Problem statement}
In the frequency regulation of optimal load control problem, only the controllable loads $d=\mathrm{col}\{d_1,...,d_n\}$ are managed to offset the change of uncontrollable power injection. There are mainly two control goals:

1. The synchronized frequency of the power network should be restored to nominal value;

2. The controllable loads should be adjusted in an economical manner without violating the operational constraints for safety, including power limits, line thermal limits and active power flow balance at each bus.

From the perspective of physical phenomenon, nonsmoothness actually emerges in many cases in smart grid, such as cost of energy storage systems, tiered price in energy market and cost of generators considering valve point loading effect \cite{tieredprice,bess2,valve,air}. For example, when aiming to control air-conditioner loads with minimum disruption, the cost function can be modelled as $f_i(d_i)=a_id_i^2+b_i\|d_i+c_i\|_1$ where $a_i$, $b_i$ and $c_i$ are coefficients \cite{air}.
Therefore, in the control goal 2, the regulation cost for each load is a convex but nonsmooth function, which can be split into three terms as $f_i(d_i)=f_{i,0}(d_i)+f_{i,1}(d_i)+f_{i,2}(d_i)$ to facilitate our load controller design. To be specific, $f_{i,0}$ is in a widely adopted form, e.g., quadratic function. $f_{i,1}$ is an indicator function that incorporates load power limits, and defined as
\begin{equation*}
   f_{i,1}(d_{i})=\begin{cases}
     0, & \mbox{if }  d_{i} \in \Omega_i\\
     +\infty, & \mbox{otherwise}
   \end{cases}
\end{equation*}
where $\Omega_i$ can be expressed as $\underline{d}_i\leq d_{i}\leq \overline{d}_i$ with $\Omega=\bigcap_{i=1}^n{\Omega_i}$.
$f_{i,2}$ is separated from $f_i$ to describe the other nonsmoothness apart from indicator function $f_{i,1}$, e.g., $l_1$ norm. Clearly, both of $f_{i,1}$ and $f_{i,2}$ are lower semicontinuous convex.

Aiming at the above two goals, the nonsmooth OLC problem is formulated as follows
\begin{subequations}\label{olc}
  \begin{alignat}{4}
\hspace{-6pt}\min\limits_{d,\theta,\omega}\quad\,\,\; &\sum_{i\in\mathcal{N}}f_i(d_i) \label{olca} \\
 \hspace{-6pt}   \text{s.t.} \quad {d_i} &{=P_i^{in}\hspace{-1.5pt}-\hspace{-1.5pt}D_i\omega_i+\sum_{{k:ki\in\mathcal{E}}}\hspace{-4pt}B_{ki}\theta_{ki}-\hspace{-2pt}\sum_{{j:ij\in\mathcal{E}}}\hspace{-5pt}B_{ij}\theta_{ij},\,\,\forall i\in \mathcal{N}} \label{olcb}\\
\hspace{-6pt}   \qquad {\omega_i}&{=0,}\qquad {\forall i\in\mathcal{N}}\label{olchoujia}\\
\hspace{-6pt}   \qquad\underline{P}_{ij}&\leq B_{ij}\theta_{ij}\leq \overline{P}_{ij}, \qquad \forall {ij\in\mathcal{E}} \label{olcc}
  \end{alignat}
\end{subequations}
where $\theta=\mathrm{col}\{\theta_1,...,\theta_n\}$ and $\omega=\mathrm{col}\{\omega_1,...,\omega_n\}$. $\underline{P}_{ij}$ and $\overline{P}_{ij}$ represent the lower and upper line thermal limits at line $ij$, respectively. In this formulation, \eqref{olcb} ensures the local power balance at each bus such that the control goal 1 is achieved. Clearly, the global power balance is guaranteed under \eqref{olcb} if summing it up over all the buses.

In the sequel, some assumptions are listed for the wellposedness of the nonsmooth OLC problem.

\textbf{\emph{Assumption 1}} (Convexity): $f_{i,0}(d_i)$ is $\beta$-strongly convex with $\beta>1$ for all $i\in\mathcal{N}$.

{\textbf{\emph{Assumption 2}} (Slater's condition): For all $i\in\mathcal{N}$, there exists $\tilde{d}=\mathrm{col}\{\tilde{d}_1,...,\tilde{d}_n\}\in relint(\Omega)$ and $\tilde{\theta}=\mathrm{col}\{\tilde{\theta}_1,...,\tilde{\theta}_n\}\in\mathbb{R}^n$ such that $0=P_i^{in}-\tilde{d}_i+\sum_{k:ki\in\mathcal{E}}B_{ki}\tilde{\theta}_{ki}-\sum_{j:ij\in\mathcal{E}}B_{ij}\tilde{\theta}_{ij}$ and $\underline{P}_{ij}< B_{ij}\tilde{\theta}_{ij}< \overline{P}_{ij}$.}

{As a consequence of Assumptions 1 and 2, the strong duality holds for OLC problem in \eqref{olc} and an optimal primal-dual pair exists. That is, the solution set of problem \eqref{olc} is nonempty.}

\begin{remark}
It is common to model the dispatching cost, i.e., $f_{i,0}(d_i)$, in the quadratic form, which is strongly convex \cite{hill,quad,ercixing2}. Yet it seems that Assumption 1 is still a bit restrictive. Actually, for those problems whose strong convex coefficient of the objective function is $\beta\in (0,1]$, {Assumption 1} can still be guaranteed by transforming the primal problem, e.g., $\min f(x)$ can be equivalently transformed into $\min kf(x)$, where $k>\frac{1}{\beta}$.
\end{remark}

\section{Algorithm Design}
It should be pointed out that one of the main challenges in this work comes from solving the nonsmooth OLC problem, in particular, the nonsmoothness from the cost function in terms of $f_{i,1}(d_i)$ and $f_{i,2}(d_i)$. Since the objective function can be separated into smooth and nonsmooth parts, our main idea is to make full use of the smooth part $f_{i,0}$, and to avoid directly calculating the subgradient of the complicated nonsmooth parts $f_{i,1}+f_{i,2}$ . Before the presentation of DPPD algorithm design, the OLC problem is reformulated for the interpretation of power system dynamics as part of optimization algorithm.
\subsection{Reformulation}
 It is observed from problem \eqref{olc} that constraints \eqref{olcb} and \eqref{olchoujia} impose both power balance requirement on each bus and the demand of frequency restoration to zero, which exactly belong to primary and secondary frequency control, respectively. In addition, the objective function \eqref{olca} and constraint \eqref{olcc} are concerning with tertiary frequency control, namely ED. Therefore, problem \eqref{olc} forms a three-tier control and optimization problem on different timescales. However, only the load demand $d$ can be adjusted, while $\theta$ reacts to the changes caused by $d$ according to \eqref{dynamics}. Therewith, it poses a requirement that the changes of $d$ in the designed optimization algorithm will not lead to the instability of power system dynamics \eqref{dynamics} which may be operated on a different timescale. Then, for the purpose of connecting the power system dynamics \eqref{dynamics} with optimization algorithm, we introduce an auxiliary variable $\hat{\theta}=\mathrm{col}\{\hat{\theta}_1,...,\hat{\theta}_n\}$ named virtual phase angle to replace $\theta$, and substitute $B_{ij}\theta_{ij}$ with $P_{ij}$ in \eqref{olcb}. Following with \cite{zch1,vfolc,pd_caid,chenxin,hill,zhaojian2}, the nonsmooth OLC problem \eqref{olc} is reformulated as
\begin{subequations}\label{rolc}
  \begin{alignat}{3}
   \hspace{-8pt}&\min\limits_{d,\hat{\theta},\omega,P}\quad \sum_{i\in\mathcal{N}}f_i(d_i)+\frac{D_i}{2}\omega_i^2 \label{rolca} \\
    \hspace{-8pt}& \text{s.t.} \;\;\; 0=P_i^{in}-D_i\omega_i-d_i+\hspace{-2pt}\sum_{{k:ki\in\mathcal{E}}}\hspace{-6pt}P_{ki}-\hspace{-2pt}\sum_{{j:ij\in\mathcal{E}}}\hspace{-6pt}P_{ij}, \forall i\in \mathcal{N}  \label{rolcb}\\
    \hspace{-8pt}&\qquad 0=P_i^{in}-d_i+\sum_{{k:ki\in\mathcal{E}}}\hspace{-6pt}B_{ki}\hat{\theta}_{ki}-\sum_{{j:ij\in\mathcal{E}}}\hspace{-6pt}B_{ij}\hat{\theta}_{ij},\forall i\in \mathcal{N} \label{rolcc}\\
   &\qquad\underline{P}_{ij}\leq B_{ij}\hat{\theta}_{ij}\leq \overline{P}_{ij}, \qquad \forall {ij\in\mathcal{E}} \label{rolcd}
  \end{alignat}
\end{subequations}
where $P=\mathrm{col}\{P_{ij}\}_{ij\in\mathcal{E}}$. $\frac{D_i\omega_i^2}{2}$ is penalized into the objection function to restore the frequency. \eqref{rolcb} is introduced to guarantee that the difference of virtual phase angle $\hat{\theta}_{ij}$ is equivalent to the true one, i.e., $\theta_{ij}$ for all $ij\in\mathcal{E}$ at steady state. Besides, the introduced $\frac{D_i\omega_i^2}{2}$ and \eqref{rolcb} are the keys to link the power system dynamics \eqref{dynamics} with the optimization algorithm. For the sake of simplicity, problem \eqref{rolc} is referred to as the reformulated optimal load control (ROLC) problem.

For problem \eqref{rolc}, the Lagrangian function is given by
\begin{eqnarray}\label{LG}
\hspace{-15pt}&&\hspace{-8pt}L(d,\hat{\theta},\omega,P,\lambda,\mu,\nu^-,\nu^+)=\sum_{i\in\mathcal{N}}f_i(d_i)+\frac{D_i}{2}\omega_i^2  \nonumber\\
  \hspace{-15pt}&&\hspace{-8pt}\quad+\sum\limits_{i\in\mathcal{N}}\lambda_i(P_i^{in}-D_i\omega_i-d_i+\sum\limits_{{k:ki\in\mathcal{E}}}P_{ki}-\sum\limits_{{j:ij\in\mathcal{E}}}P_{ij}) \nonumber \\
  \hspace{-15pt}&&\hspace{-8pt}\quad+\sum\limits_{i\in\mathcal{N}}\mu_i(P_i^{in}-d_i+\sum\limits_{{k:ki\in\mathcal{E}}}B_{ki}\hat{\theta}_{ki}-\sum\limits_{{j:ij\in\mathcal{E}}}B_{ij}\hat{\theta}_{ij} )\nonumber  \\
  \hspace{-15pt}&&\hspace{-8pt}\quad+\sum\limits_{{ij\in\mathcal{E}}}\nu_{ij}^-(\underline{P}_{ij}-B_{ij}\hat{\theta}_{ij})+\sum\limits_{{ij\in\mathcal{E}}}\nu_{ij}^+(B_{ij}\hat{\theta}_{ij}-\overline{P}_{ij}),
\end{eqnarray}
where $\lambda_i$, $\mu_i$, $\nu_{ij}^-$ and $\nu_{ij}^+$ are Lagrangian multipliers associated with constraints \eqref{rolcb}, \eqref{rolcc} and \eqref{rolcd}, respectively. $\lambda=\mathrm{col}\{\lambda_1,...,\lambda_n\}$, $\mu=\mathrm{col}\{\mu_1,...,\mu_n\}$, $\nu^-=\mathrm{col}\{\nu_{ij}^-\}_{{ij}\in\mathcal{E}}$ and $\nu^+=\mathrm{col}\{\nu_{ij}^+\}_{{ij}\in\mathcal{E}}$.

Usually, the primal-dual method is employed in the algorithm design to solve the above Lagrangian function. However, it is not a good choice to directly use subgradient to design the partial primal-dual algorithm, due to the nonsmoothness in \eqref{LG}. Specifically, on one hand, the employment of discontinuous subgradient $\partial f(d)$ makes it inapplicable for most subgradient-based algorithms to solve nonsmooth OLC problem, because the existence of solution is not guaranteed in the Carath\'{e}odory sense and violent vibration may happen during transients. On the other hand, the convergence rate of subgradient-based algorithms is relatively slow, and they do not leverage the smooth+nonsmooth structure in the algorithm design either. As for the other technique based on proximal approach, it mainly confronts with the un-proximal characteristic (inefficiently computable) if the nonsmooth function is composite. To address this problem, a three-operator splitting (TOS) method with its variants is developed in \cite{davis2017three,yuewei,randprox}. The core of TOS is further splitting the nonsmooth part into two separated parts, e.g., $f_{i,1}$ and $f_{i,2}$, and introducing an auxiliary variable $\eta_i$ to track the gradient $\partial f_{i,2}$. Then, part of complicated nonsmooth cost function, i.e., $f_{i,2}$, is separated out and its subgradient is estimated by the introduced variable. While the remained single smooth plus single nonsmooth parts, i.e., $f_{i,0}+f_{i,1}$ are dealt with proximal gradient descent approach.

\subsection{DPPD Algorithm}

As has been mentioned in the above, the main idea of DPPD is to leverage the separable property of objective function and construct a continuous estimation of partial subgradient. Subsequently, the optimization algorithm follows from a modified primal-dual framework where the whole information projections associated with dual variables are used in the dynamics of primal variables.
In the following, we will briefly explain the design procedure of TOS in continuous-time version, then propose DPPD algorithm for the ROLC problem \eqref{rolc}.

Suppose there exists an optimal auxiliary variable $\eta_i^*$ such that $-\kappa \eta_i^*$ is a subgradient in $\partial f_{i,2}(d_i^*)$, where $\kappa$ is a constant stepsize and $d_i^*$ is an optimal solution. Then, one has the following two differential inclusions at a saddle point of Lagrangian function \eqref{LG}, i.e.,
\begin{align*}
     d_i^{*}-\kappa \eta_i^*  &\in d_i^{*}+\partial f_{i,2}(d_i^{*}),\\
      d_i^{*}-\nabla f_{i,0}({d_{i}^{*}})+\kappa\eta_i^*+\lambda_i^*+\mu_i^* &\in d_i^{*}+\partial f_{i,1}(d_i^{*}),
\end{align*}
where $\lambda_i^*$ and $\mu_i^*$ are optimal dual variables of Lagrangian function \eqref{LG}. By viewing $I+\partial f_{i,1}$ and $I+\partial f_{i,2}$ on the right hand side of differential inclusions as two operators, and in virtue of the property of proximal operator shown in Section \uppercase\expandafter{\romannumeral2}-B, i.e., $\mathbf{prox}_{\tau f}(\cdot)=(I+\tau \partial f)^{-1}(\cdot)$, it yields
\begin{subequations}
\begin{align}
        &\hspace{-9pt}d_i^{*}=\mathbf{prox}_{f_{i,2}}(d_i^{*}-\kappa \eta_i^*),\label{splita}\\
        &\hspace{-9pt}d_i^{*}=\mathbf{prox}_{f_{i,1}}(d_i^{*}-\nabla f_{i,0}(d_i^{*})+\kappa\eta_i^*+\lambda_i^*+\mu_i^*).\label{splitb}
\end{align}
\end{subequations}

Accordingly, in light of \eqref{splita}, the dynamics for subgradient estimation of $\partial f_{i,2}(d_i)$ is constructed as
$$
   \dot{\eta}_i=\rho_{\eta_i}\Big[\mathbf{prox}_{f_{i,2}}(d_i-\kappa\eta_i)-d_i\Big],
$$
where $\rho_{\eta_i}$ is a constant stepsize. Next, the partial subgradient estimation $\eta_i$ replaces the original subgradient $\partial f_{i,2}(d_i)$ and helps to seek the equilibrium point $d_i^{*}$ in \eqref{splitb}. Then, the dynamics of load control $d_i$ is designed as
\begin{equation*}
  \begin{aligned}
  \dot{d}_i=\rho_{d_{i}} \Big[\mathbf{prox}_{f_{i,1}}\Big(&d_i-\nabla f_{i,0}(d_i)+\kappa \eta_i +\lambda_i+u_{i,1}\\
  &+\mu_i+u_{i,2}\Big)-d_i\Big],
  \end{aligned}
\end{equation*}
where $\rho_{d_i}$ is a constant stepsize, $u_{i,1}=P_i^{in}-d_i-D_i\omega_i+\sum\nolimits_{{k:ki\in\mathcal{E}}}P_{ki}-\sum\nolimits_{{j:ij\in\mathcal{E}}}P_{ij}$ and
$u_{i,2}=P_i^{in}-d_i+\sum\nolimits_{{k:ki\in\mathcal{E}}}B_{ki}\hat{\theta}_{ki}-\sum\nolimits_{{j:ij\in\mathcal{E}}}B_{ij}\hat{\theta}_{ij}$. It is worth noting that, by viewing $\lambda_i+u_{i,1}$ and $\mu_i+u_{i,2}$ as projections on the set of real numbers, i.e., $\mathcal{P}_{\mathbb{R}}(\lambda_i+u_{i,1})$ and $\mathcal{P}_{\mathbb{R}}(\mu_i+u_{i,2})$, then $u_{i,1}$ and $u_{i,2}$ can be considered as additional control inputs\textcolor{blue}{\footnote{\textcolor{blue}{Usually, the real-time power injection $P_i^{in}$ used in \eqref{dppd} is unknown in advance. In \cite{chenxin,zhaojian2}, they propose methods to avoid the necessity of knowing $P_i^{in}$, relying on the estimation of $P_i^{in}$. Please refer to \cite{chenxin} and \cite{zhaojian2} for details.}}} which are introduced due to the employment of the whole information projections associated with dual variables in the dynamics of primal variable $d_i$. While prior works in \cite{zch1,vfolc,zhaojian2}, that follow from the classical primal-dual framework, use dual variables solely to derive the dynamics of primal variables (e.g., $\dot{d}_i$). This is another difference in our work in terms of the construction of algorithm.

Since \eqref{splitb} can be equivalently written by projection as
\begin{equation*}
d_i^{*}=\mathcal{P}_{\Omega_i}(d_i^{*}-\nabla f_{i,0}(d_i^{*})+\kappa\eta_i^*+\lambda_i^*+\mu_i^*),
\end{equation*}
it indicates that projected gradient formulation can be viewed as a special case of \eqref{splitb}. Therefore, it is a more general perspective to study nonsmooth OLC problem with proximal approach. It shall be also pointed out that if the split nonsmooth part $f_{i,2}$ in the objective function is still un-proximal, one can introduce multiple auxiliary variables as $\eta_i$ does in \eqref{splita} to deal with it.

 \begin{figure*}[t]
   \centering
   \includegraphics[scale=0.55]{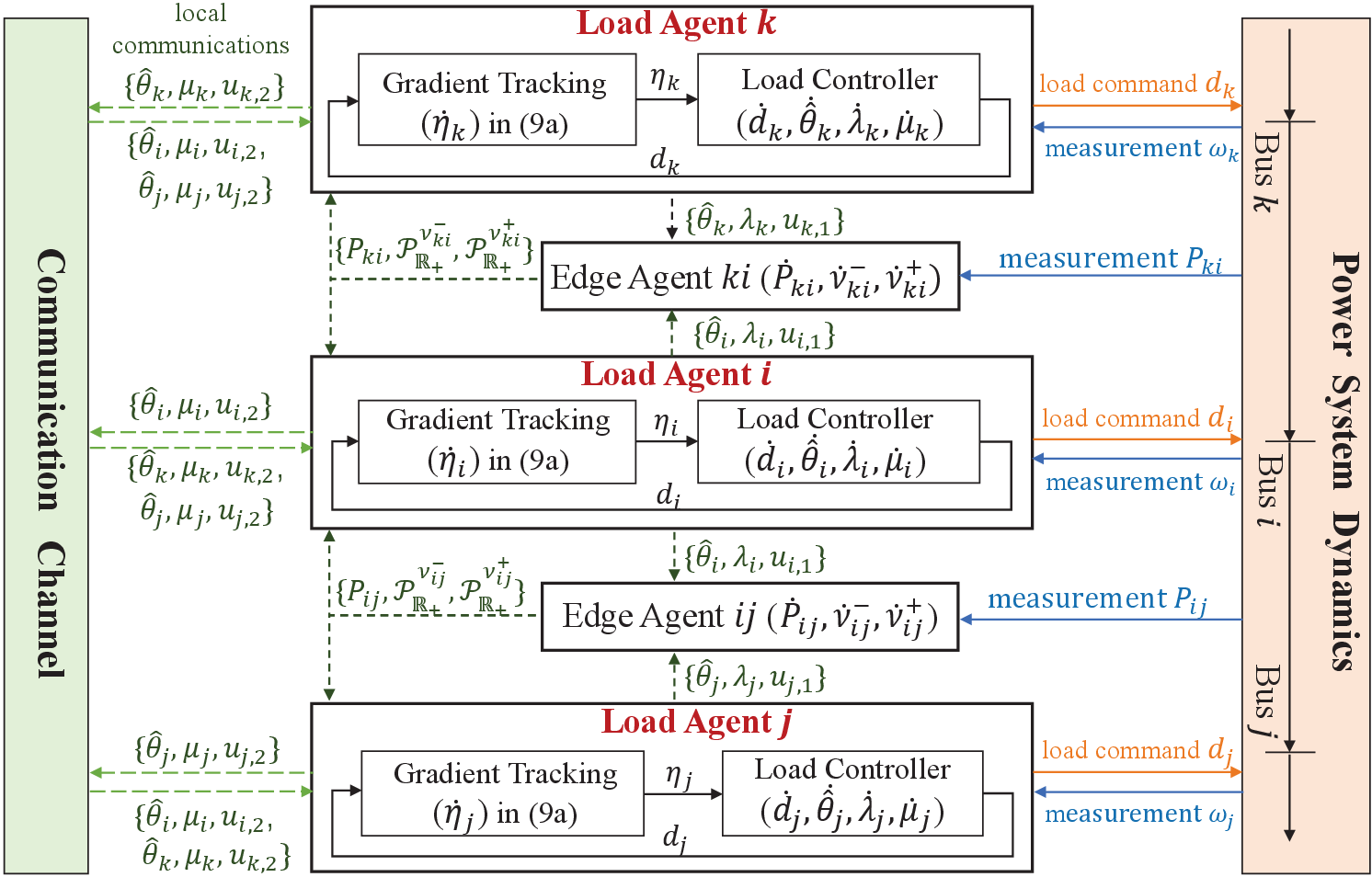}
   \caption{The architecture of DPPD algorithm}\label{architecture}
 \end{figure*}

Then, based on the above preparations, the overall solution algorithm named DPPD is proposed based on a modified primal-dual framework as follows:
 \begin{subequations}\label{dppd}
   \begin{eqnarray}
       \dot{\eta}_i\hspace{-8pt}&=&\hspace{-8pt}\rho_{\eta_i}\Big[\mathbf{prox}_{f_{i,2}}(d_i-\kappa\eta_i)-d_i\Big],  \label{dppda}\\
       \dot{d}_i\hspace{-8pt}&=&\hspace{-8pt}\rho_{d_{i}} \Big[\mathbf{prox}_{f_{i,1}}\Big(d_i-\nabla f_{i,0}(d_i)+\kappa \eta_i +\lambda_i+u_{i,1}\nonumber\\
       \hspace{-8pt}&&\qquad\qquad\quad+\mu_i+u_{i,2}\Big)-d_i\Big] ,      \label{dppdb}\\
        {\omega}_i\hspace{-8pt}&=&\hspace{-8pt}\lambda_i,\label{dppdc}\\
        \dot{\hat{\theta}}_i\hspace{-8pt}&=&\hspace{-8pt}\rho_{\hat{\theta}_i}\Big[\sum\nolimits_{{ij\in\mathcal{E}}}B_{ij}({\mathcal{P}_{\mathbb{R}_+}^{\nu_{ij}^-}}-{\mathcal{P}_{\mathbb{R}_+}^{\nu_{ij}^+}})\nonumber\\
        \hspace{-8pt}&&\quad+\sum\nolimits_{{ki\in\mathcal{E}}}B_{ki}({\mathcal{P}_{\mathbb{R}_+}^{\nu_{ki}^+}}-{\mathcal{P}_{\mathbb{R}_+}^{\nu_{ki}^-}})\nonumber\\
        \hspace{-8pt}&&\quad+\sum\nolimits_{{j:ij\in\mathcal{E}}}B_{ij}(\mu_i+u_{i,2}-\mu_j-u_{j,2})\nonumber\\
        \hspace{-8pt}&&\quad-\sum\nolimits_{{k:ki\in\mathcal{E}}}B_{ki}(\mu_k+u_{k,2}-\mu_i-u_{i,2})\Big],\label{dppdd}\\
       \dot{P}_{ij}\hspace{-8pt}&=&\hspace{-8pt}\rho_{P_{ij}}(\lambda_i+u_{i,1}-\lambda_j-u_{j,1}),\label{dppde}\\
         \dot{\lambda}_i\hspace{-8pt}&=&\hspace{-8pt}\rho_{\lambda_{i}}u_{i,1},\label{dppdf}\\
        \dot{\mu}_i\hspace{-8pt}&=&\hspace{-8pt}\rho_{\mu_i}u_{i,2},\label{dppdg}\\
        \dot{\nu}_{ij}^-\hspace{-8pt}&=&\hspace{-8pt}\rho_{{\nu}_{ij}^-}\Big[{\mathcal{P}_{\mathbb{R}_+}(\nu_{ij}^-+\underline{P}_{ij}-B_{ij}\hat{\theta}_{ij})}-\nu_{ij}^-\Big],\label{dppdh}\\
        \dot{\nu}_{ij}^+\hspace{-8pt}&=&\hspace{-8pt}\rho_{{\nu}_{ij}^+}\Big[{\mathcal{P}_{\mathbb{R}_+}(\nu_{ij}^++B_{ij}\hat{\theta}_{ij}-\overline{P}_{ij})}-\nu_{ij}^+\Big],\label{dppdi}
   \end{eqnarray}
 \end{subequations}
where $\rho_{\{\cdot\}}$ are stepsizes for each dynamics, and for brevity,
$\mathcal{P}_{\mathbb{R}_+}^{\nu_{ij}^-}:= \mathcal{P}_{\mathbb{R}_+}(\nu_{ij}^-+\underline{P}_{ij}-B_{ij}\hat{\theta}_{ij})$ and
$\mathcal{P}_{\mathbb{R}_+}^{\nu_{ij}^+}:=\mathcal{P}_{\mathbb{R}_+}(\nu_{ij}^++B_{ij}\hat{\theta}_{ij}-\overline{P}_{ij})$.

\begin{remark}
    If \eqref{dppdf} and \eqref{dppdg} are equivalently rewritten as follows: $$\dot{\lambda}_i=\rho_{\lambda_i}[\mathcal{P}_{\mathbb{R}}(\lambda_i+u_{i,1})-\lambda_i],$$ $$ \dot{\mu}_i=\rho_{\mu_i}[\mathcal{P}_{\mathbb{R}}(\mu_i+u_{i,2})-\mu_i],$$
     one can observe from \eqref{dppdb}, \eqref{dppdd} and \eqref{dppde} that the whole information projections, i.e., $\mathcal{P}_{\mathbb{R}}(\lambda_i+u_{i,1})$ and $\mathcal{P}_{\mathbb{R}}(\mu_i+u_{i,2})$, are utilized in the gradient descent process of primal variables, rather than solely relying on the dual variables $\lambda_i$ and $\mu_i$. Likewise, $\mathcal{P}_{\mathbb{R}_+}^{\nu_{ij}^-}$ and $\mathcal{P}_{\mathbb{R}_+}^{\nu_{ij}^+}$ are employed in \eqref{dppdd} instead of only using $\nu_{ij}^-$ and $\nu_{ij}^+$.
     This is why we claim that DPPD is derived from a modified primal-dual framework.

     Moreover, by substituting $\dot{\lambda}_i=\rho_{\lambda_i}u_{i,1}$ and $\dot{\mu}_i=\rho_{\mu_i}u_{i,2}$ into \eqref{dppdb}, it yields a PD-type protocol with regard to dual variables as follows:
\begin{equation*}
\begin{aligned}
  \dot{d}_i=&\rho_{d_{i}} \Big[\mathbf{prox}_{f_{i,1}}\Big(d_i-\nabla f_{i,0}(d_i)+\kappa \eta_i +\lambda_i+\frac{1}{\rho_{\lambda_{i}}} \dot{\lambda}_i
      \\
      &\qquad\qquad\qquad+\mu_i+\frac{1}{\rho_{\mu_i}}\dot{\mu}_i  \Big)-d_i\Big].
\end{aligned}
\end{equation*}
Conversely, it yields a PI-type protocol with regard to $u_{i,1}$ and $u_{i,2}$, i.e.,
\begin{equation*}
\begin{aligned}
  \dot{d}_i=&\rho_{d_{i}} \Big[\mathbf{prox}_{f_{i,1}}\Big(d_i-\nabla f_{i,0}(d_i)+\kappa \eta_i +\rho_{\lambda_i}\int u_{i,1}+u_{i,1}\\
  &\qquad\qquad\qquad
       +\rho_{\mu_i}\int u_{i,2}+u_{i,2}\Big)-d_i\Big].
\end{aligned}
\end{equation*}

Such a PD-type or PI-type protocol also presents in dynamics of other primal optimization variables, e.g., $\dot{\hat{\theta}}_i$ and $\dot{P}_{ij}$. This is an interesting change in the protocol construction, resulting from exploiting the whole information projection associated with the dual variable in the gradient dynamics of the primal variable. Intuitively, the extra degree of freedom offered by the derivative or integral feedback is expected to bring an added benefit, e.g., enlarging the delay consensus margin \cite{TAC_madan_delay}. It is, thus, desirable to assess the improvements by such protocols in our future work.
\end{remark}

\begin{remark}
   The proposed DPPD algorithm can be interpreted as follows:
  \eqref{dppda} is designed as a feedback controller to render \eqref{splita}. \eqref{dppdb} is designed by proximal gradient descent method with the whole information projections, i.e., $\mathcal{P}_{\mathbb{R}}(\lambda_i+u_{i,1})$ and $\mathcal{P}_{\mathbb{R}}(\mu_i+u_{i,2})$, and consequently rendering \eqref{splitb}. \eqref{dppdc} is obtained by firstly minimizing the Lagrangian function \eqref{LG} over $\omega_i$, and meanwhile we use a trick resulted from ${\partial L}/{\partial \omega_i}=0$. That is, it has $\omega_i^*=\lambda_i^*$ at the equilibrium which means $\omega_i$ and $\lambda_i$ are interchangeable in the algorithm design if the proposed algorithm is a contractive operator. Then, it yields \eqref{dppdc}.
  \eqref{dppdd}-\eqref{dppdi} are procured along with the optimization of Lagrangian function \eqref{LG} with the aforementioned modified primal-dual dynamics. At this moment, considering the trick that $\omega_i$ and $\lambda_i$ are interchangeable, the optimization algorithm (specifically \eqref{dppde} and \eqref{dppdf}) is related to the dynamics of the power network as in \eqref{3b} and \eqref{3c} on the same timescale but with a PD-type protocol, if the stepsizes $\rho_{P_{ij}}=B_{ij}$ and $\rho_{\lambda_i}=\frac{1}{M_i}$. Additionally, note that the loads considered here are frequency insensitive, which means $M_i=0$ for $i \in \mathcal{N_L}$. Thus, if $\rho_{\lambda_i}\rightarrow +\infty$, we have $0=P_i^{in}-d_i-D_i\omega_i+\sum_{{k:ki\in\mathcal{E}}}P_{ki}-\sum_{{j:ij\in\mathcal{E}}}P_{ij}$ from \eqref{dppdf} for $i\in\mathcal{N_L}$, which is exactly the dynamics of load in \eqref{3d}.
\end{remark}

\begin{remark}
  {Since the smoothness of DPPD algorithm is ensured, and consequently it avoids the nonconvexity caused by the projection of discontinuous subgradient. Hence, the existence of a Carath\'{e}odory solution to DPPD algorithm is guaranteed.}
\end{remark}
\section{Optimality and Convergence Analysis}
In this section, an equilibrium point of the proposed DPPD algorithm is proved to be an optimal solution of ROLC problem that satisfies the Karush-Kuhn-Tucker (KKT) condition. Later, the convergence result is analyzed via Lyapunov stability and LaSalle invariance principle. For simplicity, the stepsizes $\rho_{ \{\cdot\} }$ associated with the dynamics \eqref{dppda}-\eqref{dppde} are set to be 1, while those in \eqref{dppdf}-\eqref{dppdi} are set to be $\frac{1}{2}$, and $\kappa\in(0,1)$.

\subsection{Optimality}
Let $\eta=\mathrm{col}\{\eta_1,...,\eta_n\}$ (doing likewise with $u_1$, $u_2$ and $P^{in}$), $D=\mathrm{diag}\{D_1,...,D_n\}$ and $B=\mathrm{diag}\{B_{ij}\}_{ij\in\mathcal{E}}$. Denote $(\eta^*,{d^{*}},{\hat{\theta}}^*,\omega^*,P^*,\mu^*,\nu^{-*},\nu^{+*})$ as an equilibrium point of DPPD algorithm. Firstly, a modified optimality condition is given based on the Karush-Kuhn-Tucker (KKT) condition.
\par \emph{Lemma 1}:
    $({d^{*}},\hat{\theta}^*,\omega^*,P^*)$ is an optimal solution to the ROLC problem if and only if $(\mu^*,\nu^{-*},\nu^{+*})$ exists such that
    \begin{subequations}\label{kkt}
        \begin{alignat}{8}
            \mathbf{0}\in&\hspace{2pt}\partial_{d}(\sum_{i=1}^{n}f_i(d_i^{*}))-{\lambda}^*-{\mu}^*,\label{kkta}\\
            \omega^*=&\hspace{2pt}\lambda^*,\label{kktc}\\
            \mathbf{0}=&\hspace{2pt}CB(\nu^{-*}-\nu^{+*})+CBC^T\mu^*,\label{kktb}\\
            \mathbf{0}=&\hspace{2pt}C^T\lambda^*,\label{kktd}\\
            \mathbf{0}=&\hspace{2pt}P^{in}-{d^{*}}-D\omega^*-CP^*,\label{kkte}\\
            \mathbf{0}=&\hspace{2pt}P^{in}-{d^{*}}-CBC^T\hat{\theta}^*,\label{kktf}\\
            \underline{P}-BC^T\hat{\theta}^*\in&\hspace{2pt}\mathcal{C}_{\mathbb{R}_{+}^l}(\nu^{-*}),\label{kktg}\\
            BC^T\hat{\theta}^*-\overline{P}\in&\hspace{2pt}\mathcal{C}_{\mathbb{R}_{+}^l}(\nu^{+*}).\label{kkth}
        \end{alignat}
    \end{subequations}

\begin{IEEEproof}
  Lemma 1 can be obtained from the KKT condition with the definition of normal cone. Specifically, \eqref{kkta}-\eqref{kktd} are stationarity conditions along with the minimization of Lagrangian function \eqref{LG}. \eqref{kkte} and \eqref{kktf} are the primal feasibility of auxiliary constraints \eqref{rolcb} and \eqref{rolcc}, respectively. \eqref{kktg} and \eqref{kkth} denote the complementary slackness and dual feasibility of active power flow constraint \eqref{rolcd}. This ends the proof.
\end{IEEEproof}

 According to the equilibrium state of \eqref{dppd}, the associated compact form can be written as
\begin{subequations}\label{sstate}
  \begin{eqnarray}
    -\kappa\eta^*\hspace{-6pt} &\in&\hspace{-6pt} \partial f_2(d^{*}), \label{sstatea}\\
    \mathbf{0}\hspace{-6pt}&\in&\hspace{-6pt}\nabla f_0(d^{*})+\partial f_1(d^{*})-\kappa\eta^*-\lambda^*-\mu^*,\label{sstateb}\\
    \omega^*\hspace{-6pt} &=& \hspace{-6pt}\lambda^*, \label{sstated}\\\mathbf{0}\hspace{-6pt} &=& \hspace{-6pt}CB(\nu^{-*}-\nu^{+*})+CBC^T\mu^*, \label{sstatec}\\
    \mathbf{0}\hspace{-6pt} &=& \hspace{-6pt}C^T\lambda^*, \label{sstatee}\\
    \mathbf{0}\hspace{-6pt} &=& \hspace{-6pt}P^{in}-d^{*}-D\omega^*-CP^*, \label{sstatef}\\
    \mathbf{0}\hspace{-6pt} &=& \hspace{-6pt}P^{in}-d^{*}-CBC^T\hat{\theta}^*, \label{sstateg}\\
  \mathbf{0}\hspace{-6pt}&\in&\hspace{-6pt}\mathcal{C}_{\mathbb{R}_{+}^l}(\nu^{-*})-\underline{P}+BC^T\hat{\theta}^*,\label{sstateh}\\
            \mathbf{0}\hspace{-6pt}&\in&\hspace{-6pt}\mathcal{C}_{\mathbb{R}_{+}^l}(\nu^{+*})-BC^T\hat{\theta}^*+\overline{P},\label{sstatei}
  \end{eqnarray}
\end{subequations}
where $\nabla f_0(d^{*})=\{\nabla f_{1,0}(d_1^{*}),...,\nabla f_{n,0}(d_n^{*})\}$, $\partial f_j(d^{*})=\{\partial f_{1,j}(d_1^{*}),...,\partial f_{n,j}(d_n^{*})\}$, $j\in\{1,2\}$.
Following with Lemma 1 and the above equilibrium state, the equivalence between the equilibrium point and the optimal solution of ROLC problem is discussed.

\begin{theorem}
  Under Assumptions 1 and 2, $({d^{*}},\hat{\theta}^*,\omega^*=\mathbf{0},P^*)$ is an optimal solution of the ROLC problem, with $\eta^*$ being the optimal gradient tracking of $\partial f_2({d^{*}})$ and $(\mu^*,\nu^{-*},\nu^{-*})$ being an optimal solution to the dual problem.
\end{theorem}

\begin{IEEEproof}
As for the sufficiency firstly, it is straightforward that \eqref{kkta} holds by substituting \eqref{sstatea} into \eqref{sstateb}. \eqref{kktb}-\eqref{kkth} hold for \eqref{sstatec}-\eqref{sstatei}, respectively. Besides, we have that following conclusions. First, in light of the property of the incidence matrix that $\mathbf{1}^T C=\mathbf{0}^T$, it has $\omega^*=\lambda^*=\epsilon_1\mathbf{1}$ from \eqref{sstated} and \eqref{sstatee} where $\epsilon_1$ is a constant. Second, by subtracting \eqref{sstatef} from \eqref{sstateg} and multiplying $\mathbf{1}^T$ on both sides, one has $\mathbf{1}^TD\omega^*=0$. Note that $D$ is a positive definite diagonal matrix and $\omega^*=\epsilon_1\mathbf{1}$, then it has $\omega^*=\mathbf{0}$. That is to say, the frequency deviations on all the buses are restored to zero. Further, by substituting $\omega^*=\mathbf{0}$ into \eqref{sstatef}, it yields $P^*=BC^T\theta^*=BC^T\hat{\theta}^*$ from \eqref{sstatef} and \eqref{sstateg}, leading to the result that $\theta^*=\hat{\theta}^*+\epsilon_2\mathbf{1}$ where $\epsilon_2$ is a constant. That means $\theta_{ij}^*=\hat{\theta}_{ij}^*$ for all ${ij\in\mathcal{E}}$.

Conversely, suppose that $({d}^*,\hat{\theta}^*,\omega^*,P^*,\mu^*,\nu^{-*},\nu^{+*})$ is an optimal solution. For a separable objective function, there must exist an $\eta^*\in \mathbb{R}^n$ with $\kappa$ such that the subgradient of partial nonsmooth function is trackable, i.e., $-\kappa\eta^*\in\partial f_2({d^{*}})$. Then, \eqref{sstatea} and \eqref{sstateb} are obtained. As for \eqref{sstatec}-\eqref{sstatei}, the necessity comes naturally.
\end{IEEEproof}

\begin{remark}
  With regard to an equilibrium point $({d^{*}},\hat{\theta}^*,\omega^*,P^*)$, the uniqueness of $d^*$ and $\omega^*$ are straightforward, because the objective function is strongly convex over $d$ and $\omega$, respectively. While $\hat{\theta}^*$ and $P^*$ are not always unique, particularly, in a meshed power network. In such situation, the incidence matrix $C$ is not full column rank, leading to many $P^*$ that satisfies \eqref{kkte}. Similarly, it has many different $\hat{\theta}^*$ satisfying \eqref{kktf}.
\end{remark}

The next lemma guarantees the equivalence of problem reformulation of \eqref{olc} into \eqref{rolc}, which is naturally derived from Lemma 4 in \cite{vfolc}.

\par {\emph{Lemma 2} (\hspace{-0.003pt}\cite{vfolc}):
$(d^*,\hat{\theta}^*,\omega^*,P^*)$ is an optimal solution of problem \eqref{rolc}. Then, $(d^*,\theta^*,\omega^*)$ is an optimal solution of problem \eqref{olc} if $C^T\theta^*=C^T\hat{\theta}^*$ and $BC^T\theta^*=P^*$.}

\begin{IEEEproof}
The proof is given in Appendix A.
\end{IEEEproof}

\subsection{Convergence}
Under the aforementioned preparations, we come to a position to state the main results. For the sake of presentation, we firstly show the convergence results that line thermal limits in \eqref{rolcd} are not considered in the subsequent two theorems, and then extend to the case where \eqref{rolcd} exists in Corollary.
\begin{theorem}
  If Assumptions 1 and 2 hold, the following statements are true under the proposed DPPD algorithm:

  1) Each $d_i$ remains within the load power limits, i.e., $d_i\in\Omega_i$, for any initially feasible value;

  2) The trajectory $(\eta,{d},{\hat{\theta}},P,\omega,\mu)$ generated by the proposed DPPD algorithm is bounded;

  3) {Furthermore, the trajectory $(\eta,{d},{\hat{\theta}},P,\omega,\mu)$ converges {globally} and {asymptotically} to an equilibrium point of DPPD algorithm, i.e., $(\eta^*,{d}^*,{\hat{\theta}}^*,P^*,\omega^*,\mu^*)$, where $({d}^*,{\hat{\theta}}^*,P^*,\omega^*)$ is an optimal solution to problem \eqref{rolc}.}
\end{theorem}
\begin{IEEEproof}
The proof is given in Appendix B.
\end{IEEEproof}

\begin{remark}
  Note that the exact values of damping coefficients $D_i$, $ i\in\mathcal{N}$, are not required in the sufficient conditions for the stability analysis of closed-loop system. Whereas in \cite{hill}, exact $D_i$ or their bounds, i.e., $D_i\in[D_i^{\min},D_i^{\max}]$, are used for the selection of parameters associated with the frequency controllers, leading to conservatism and inconvenient applicability in practice.
\end{remark}


Next, based on Theorem 2, the convergence rate result in terms of dynamics of DPPD algorithm is proved by contradiction.
\begin{theorem}
  With the same conditions used in Theorem 2, the convergence rate in the sense of dynamics of optimization variables in DPPD algorithm can be measured as follows:
$$\liminf (\|\dot{d}\|+\|\dot{\eta}\|+\|\dot{\hat{\theta}}\|+\|\dot{P}\|+\|\dot{\mu}\|+\|\dot{\omega}\|)=O(1/\sqrt{t}).$$
\end{theorem}

\begin{IEEEproof}
   The convergence rate analysis follows with the convergence result shown in \eqref{dot_Vb} when line thermal limits \eqref{rolcd} are not considered. By \eqref{dot_Vb} and specifying a non-negative constant $\hat{c}\leq 1-\frac{\alpha(1+\kappa)}{2}$, it has that
\begin{eqnarray}
  \hspace{-12pt}\dot{V}_b\hspace{-7pt}&\leq&\hspace{-7pt}-\|\dot{d}\|^2-(1-\frac{1+\kappa}{2\alpha\beta})(\nabla f_0(d)-\nabla f_0(d^{*}))^T\Delta d \nonumber\\
  \hspace{-12pt}\hspace{-7pt}&&\hspace{-7pt} -(1-\frac{\alpha(1+\kappa)}{2})\|\dot{\eta}\|^2-\|\mathcal{P}_{\Theta}^{\hat{\theta}}-\hat{\theta}\|^2-\|\mathcal{P}_{\mathbb{R}^l}^{P}-P\|^2\nonumber\\
  \hspace{-12pt}\hspace{-7pt}&&\hspace{-7pt}-\frac{1}{2}\|\mathcal{P}_{\mathbb{R}^n}^{\mu}-\mu\|^2-\frac{1}{2}\|\mathcal{P}_{\mathbb{R}^n}^{\omega}-\omega\|^2\nonumber\\
   \hspace{-12pt}\hspace{-7pt}&\leq&\hspace{-7pt}-\|\dot{d}\|^2-(1-\frac{\alpha(1+\kappa)}{2})\|\dot{\eta}\|^2-\|\dot{\hat{\theta}}\|^2-\|\dot{P}\|^2\nonumber\\
   \hspace{-12pt}\hspace{-7pt}&&\hspace{-7pt}-2\|\dot{\mu}\|^2-2\|\dot{\omega}\|^2\nonumber\\
   \hspace{-12pt}\hspace{-7pt}&\leq&\hspace{-7pt} -\hat{c}(\|\dot{d}\|^2+\|\dot{\eta}\|^2+\|\dot{\hat{\theta}}\|^2+\|\dot{P}\|^2 +\|\dot{\mu}\|^2+\|\dot{\omega}\|^2)\nonumber\\
  \hspace{-12pt}\hspace{-7pt}&\leq&\hspace{-7pt} -\frac{\hat{c}}{4}(\|\dot{d}\|+\|\dot{\eta}\|+\|\dot{\hat{\theta}}\|+\|\dot{P}\|+\|\dot{\mu}\|+\|\dot{\omega}\|)^2\label{rate}
\end{eqnarray}
where the second inequality is obtained by using the convexity of $f_0(\cdot)$ and the facts that $1-\frac{1+\kappa}{2\alpha\beta}\geq0$, $1-\frac{\alpha(1+\kappa)}{2}\geq 0$, $\dot{\hat{\theta}}=\mathcal{P}_{\Theta}^{\hat{\theta}}-\hat{\theta}$, $\dot{P}=\mathcal{P}_{\mathbb{R}^l}^{P}-P$, $\dot{\mu}=\frac{1}{2}(\mathcal{P}_{\mathbb{R}^n}^{\mu}-\mu)$ and $\dot{\omega}=\frac{1}{2}(\mathcal{P}_{\mathbb{R}^n}^{\omega}-\omega)$ as in \eqref{dppd}. By this upper bound of $\dot{V}_b$ in \eqref{rate}, we can obtain that the partial convergence rate of  $(\dot{d},\dot{\eta},\dot{\hat{\theta}},\dot{P},\dot{\mu},\dot{\omega})$ is $O(1/\sqrt{t})$. If it is not true, that means there exists a timestamp $t_0$ such that $\|\dot{d}\|+\|\dot{\eta}\|+\|\dot{\hat{\theta}}\|+\|\dot{P}\|+\|\dot{\mu}\|+\|\dot{\omega}\|\geq {c}/\sqrt{t}$ for all $t\geq t_0$ where ${c}$ is an any large positive constant. Then, combining with \eqref{rate}, it has
\begin{eqnarray}\label{2}
   V_b(t) \hspace{-6pt}&\leq& \hspace{-6pt} V_b(t_0)-\int_{t_0}^{t}\frac{\hat{c}}{4}\frac{{c}^2}{x}dx\nonumber\\
   \hspace{-6pt}&\leq& \hspace{-6pt} V_b(t_0)+\frac{\hat{c} c^2}{4}(\ln{t_0}-\ln t).
\end{eqnarray}
One can observe that $V_b(t)<0$ when $t$ is large enough, which contradicts with the result in Theorem 2, i.e., $V_b(t)\geq 0$. Then, it has $\liminf (\|\dot{d}\|+\|\dot{\eta}\|+\|\dot{\hat{\theta}}\|+\|\dot{P}\|+\|\dot{\mu}\|+\|\dot{\omega}\|)=O(1/\sqrt{t})$ by such an observation. This completes the proof.
\end{IEEEproof}

\begin{remark}
    Limited by the technical difficulties in convergence proofs caused by nonsmoothness and inequality capacity constraints in nonsmooth OLC problem, it is hard to obtain a linear convergence rate based on the currently constructed Lyapunov function candidate $V_b$ and the associated $\dot{V}_b$, though the strongly convex assumption is used for the objective function $f(d)$. We also notice that some related works, e.g., \cite{Proximal-Augmented-Lagrangian-Method-tac-2019,Exponential-Stability-of-Primal-Dual-Gradient-Dynamics-NaLi,Exponentially-Convergent-Algorithm-Design-for-Constrained-Distributed-Optimization-via-Nonsmooth-Approach}, developed exponential stability or convergence rate with (proximal) augmented Lagrangian. In our future work, we will explore the linear convergence rate of the proposed DPPD algorithm in nonsmooth OLC problem.
\end{remark}

{Before the end of this section, the above convergence results are extended to the case where line thermal limits \eqref{rolcd} are considered.
\begin{corollary}
  Under the same conditions in Theorem 2, the trajectory $(\eta,{d},{\hat{\theta}},P,\omega,\mu,\nu^{-},\nu^{+})$ globally asymptotically converges to an equilibrium point of DPPD algorithm when considering constraint \eqref{rolcd}. Meanwhile, the partial convergence rate of DPPD algorithm still holds as $\liminf (\|\dot{d}\|+\|\dot{\eta}\|+\|\dot{\hat{\theta}}\|+\|\dot{P}\|+\|\dot{\omega}\|+\|\dot{\mu}\|+\|\dot{\nu}^-\|+\|\dot{\nu}^+\|)=O(1/\sqrt{t})$.
\end{corollary}}

\begin{IEEEproof}
When considering the inequality capacity constraint as in \eqref{rolcd}, the dynamics of $\dot{\nu}^{-}$ and $\dot{\nu}^{+}$ start to work in DPPD algorithm. In essence, the roles of $\mu$, $\nu^-$ and $\nu^+$ are almost the same. Specifically, in the proof of statement 2 in Theorem 2, the dynamics of $\dot{\mu}$ is equivalently treated as the projection on the set of real numbers, i.e., $\dot{\mu}=u_2$ is rewritten as $\dot{\mu}=\mathcal{P}_{\mathbb{R}^n}(\mu+u_2)-\mu$. The latter formulation has a close connection with dynamics of $\dot{\nu}^{-}$ and $\dot{\nu}^{+}$ as in \eqref{dppdh} and \eqref{dppdi}. As a result, the cores of proof are in line with those in statement 2.

According to the gradient dynamics in DPPD algorithm, define the auxiliary function as
\begin{eqnarray*}
   \hspace{-6pt}&&\hspace{-6pt}\Psi^{'}(d,\hat{\theta},P,\omega,\mu,\nu^-,\nu^+)\\
   \hspace{-6pt}&=&\hspace{-6pt}f_0(d)+\frac{1}{2}\|\mathcal{P}_{\mathbb{R}^n}^{\omega}\|^2+\frac{1}{2}\|\mathcal{P}_{\mathbb{R}^n}^{\mu}\|^2+\frac{1}{2}\|\mathcal{P}_{\mathbb{R}_+^l}^{\nu^-}\|^2+\frac{1}{2}\|\mathcal{P}_{\mathbb{R}_+^l}^{\nu^+}\|^2,
\end{eqnarray*}
where $\mathcal{P}_{\mathbb{R}_+^l}^{\nu^-}:=\mathcal{P}_{\mathbb{R}_+^l}(\nu^-+\underline{P}-BC^T\hat{\theta})=\mathrm{col}\{\mathcal{P}_{\mathbb{R}_+}^{\nu_{ij}^-} \}_{ij\in\mathcal{E}}$ and $\mathcal{P}_{\mathbb{R}_+^l}^{\nu^+}:=\mathcal{P}_{\mathbb{R}_+^l}(\nu^++BC^T\hat{\theta}-\overline{P})=\mathrm{col}\{\mathcal{P}_{\mathbb{R}_+}^{\nu_{ij}^+}\}_{ij\in\mathcal{E}}$. Based on $\Psi^{'}$, the Lyapunov function candidate can be constructed as $$V_b^{'}=V_1^{'}+V_2+V_3^{'}+V_4,$$ where
\begin{eqnarray*}
  \hspace{-16pt}V_1^{'} \hspace{-6pt}&=& \hspace{-6pt} \Psi^{'}-\Psi^{'*}-\nabla_{d}^T{\Psi^{'*}}\Delta d-2{\omega^*}^T(I-D)\Delta\omega\\
  \hspace{-16pt}&&-2{\mu^*}^T\Delta\mu+2{\omega^*}^TC\Delta P+2{\mu^*}^TCBC^T\Delta \hat{\theta}\\
  \hspace{-16pt}&& -{\nu^{-*}}^T\Delta\nu^--{\nu^{+*}}^T\Delta\nu^+,\\
  \hspace{-16pt}V_3^{'}\hspace{-6pt}&=& \hspace{-6pt}V_{3.1}+V_{3.2}+\frac{1}{2}(\|\Delta\nu^-\|^2+\|\Delta\nu^+\|^2).
\end{eqnarray*}
By utilizing the tricks used in the proof of statement 2 to offset the related terms introduced by $\nu^-$ and $\nu^+$, specifically from \eqref{mu_1st} to \eqref{23}, one can obtain that $\dot{V}_b^{'}$ is negative semidefinite as follows
\begin{eqnarray}\label{withnu}
  \dot{V}_b^{'}\hspace{-7pt}&\leq&\hspace{-7pt} -\|\dot{d}\|^2-(1-\frac{1+\kappa}{2\alpha\beta})(\nabla f_0(d)- \nabla f_0(d^{*}))^T\Delta d \nonumber\\
 \hspace{-7pt} && \hspace{-7pt} -(1-\frac{\alpha(1+\kappa)}{2})\|\dot{\eta}\|^2-\|\mathcal{P}_{\Theta}^{\hat{\theta}}-\hat{\theta}\|^2-\|\mathcal{P}_{\mathbb{R}^l}^{P}-P\|^2\nonumber\\
 \hspace{-7pt} &&\hspace{-7pt}-\frac{1}{2}\|\mathcal{P}_{\mathbb{R}^n}^{\mu}-\mu\|^2-\frac{1}{2}\|\mathcal{P}_{\mathbb{R}^n}^{\omega}-\omega\|^2-\frac{1}{2}\|\mathcal{P}_{\mathbb{R}_+^l}^{\nu^-}-\nu^-\|^2\nonumber\\
 \hspace{-7pt} &&\hspace{-7pt}-\frac{1}{2}\|\mathcal{P}_{\mathbb{R}_+^l}^{\nu^+}-\nu^+\|^2.
\end{eqnarray}

To analyze the partial convergence rate of DPPD algorithm, it proceeds with \eqref{withnu} as follows:
\begin{eqnarray}\label{ratenu_1}
  \dot{V}_b^{'}\hspace{-7pt}&\leq &\hspace{-7pt} -\|\dot{d}\|^2-(1-\frac{\alpha(1+\kappa)}{2})\|\dot{\eta}\|^2-\|\dot{\hat{\theta}}\|^2-\|\dot{P}\|^2 \nonumber\\
  \hspace{-7pt}&&\hspace{-7pt}-2\|\dot{\mu}\|^2-2\|\dot{\omega}\|^2-2\|\dot{\nu}^-\|^2-2\|\dot{\nu}^+\|^2\nonumber\\
 \hspace{-7pt}&\leq& \hspace{-7pt} -\hat{c}(\|\dot{d}\|^2+\|\dot{\eta}\|^2+\|\dot{\hat{\theta}}\|^2+\|\dot{P}\|^2 +\|\dot{\mu}\|^2+\|\dot{\omega}\|^2\nonumber\\
 \hspace{-7pt}&& \hspace{-7pt}\quad\quad+\|\dot{\nu}^-\|^2+\|\dot{\nu}^+\|^2)\nonumber\\
 \hspace{-7pt}&\leq&\hspace{-7pt} -\frac{\hat{c}}{4}(\|\dot{d}\|+\|\dot{\eta}\|+\|\dot{\hat{\theta}}\|+\|\dot{P}\|+\|\dot{\mu}\|+\|\dot{\omega}\|\nonumber\\
 \hspace{-7pt}&&\hspace{-7pt}\quad\quad+\|\dot{\nu}^-\|+\|\dot{\nu}^+\|)^2.
\end{eqnarray}
Subsequently, as has been discussed in the proof of Theorem 3, one can obtain that
  $\liminf (\|\dot{d}\|+\|\dot{\eta}\|+\|\dot{\hat{\theta}}\|+\|\dot{P}\|+\|\dot{\omega}\|+\|\dot{\mu}\|+\|\dot{\nu}^-\|+\|\dot{\nu}^+\|)=O(1/\sqrt{t})$.
This completes the proof.

\end{IEEEproof}

\section{Case Studies}
In this section, the effectiveness of our proposed DPPD algorithm \eqref{dppd} is validated. First, the performance of DPPD algorithm is tested under step change and time-varying power injections. Then, we illustrate the robustness in two cases under parameter uncertainties.

\subsection{Test system}
\begin{figure}[t]
   \centering
   \includegraphics[scale=0.33]{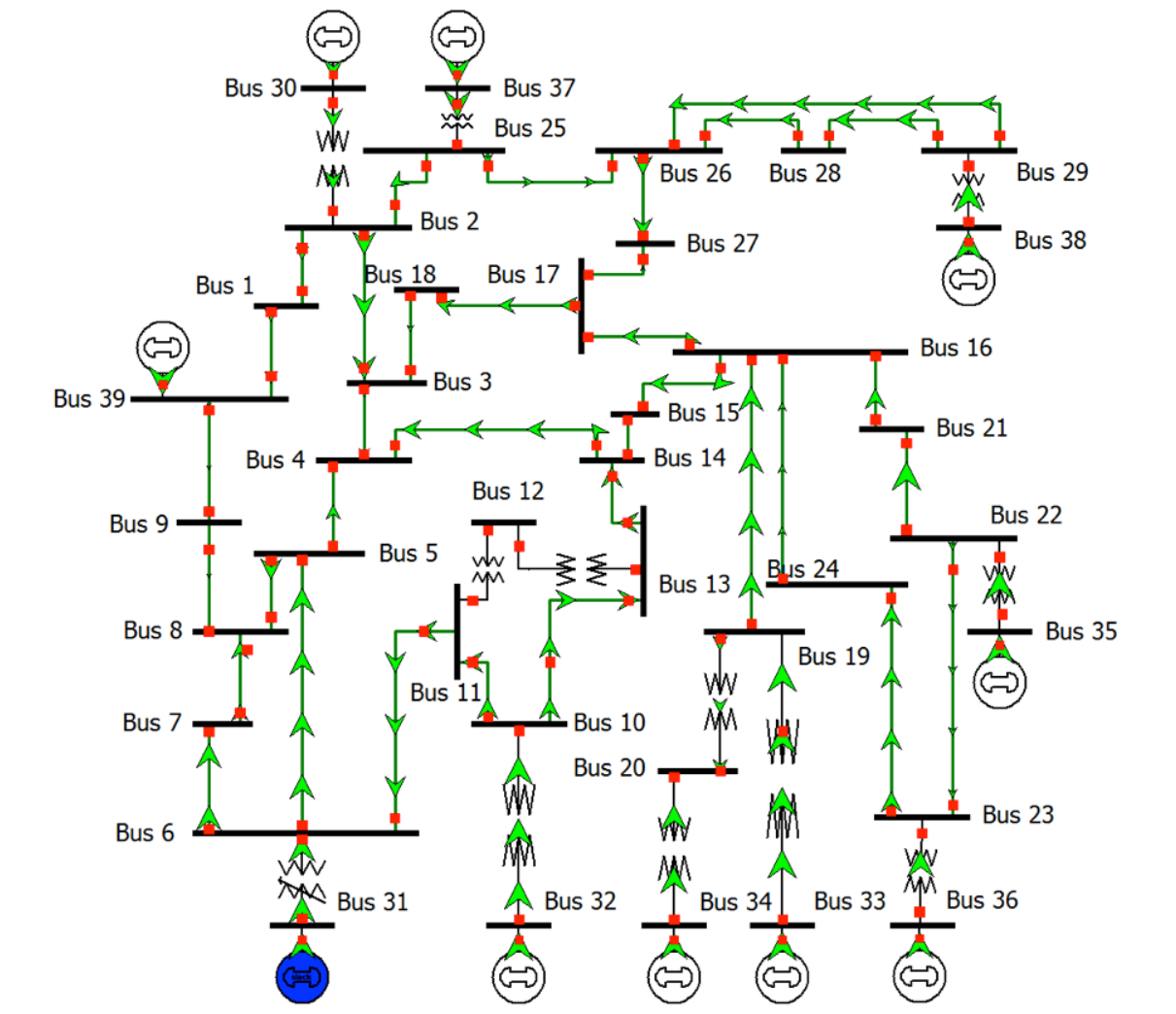}
   \caption{IEEE 39-bus power network}\label{39bus}
 \end{figure}

The simulations are carried out on Matlab R2022a with Power System Toolbox \cite{pst}. The standard IEEE 39-bus power network that consists of 10 generator buses, 29 load buses and 46 transmission lines is used as the test system shown in Fig. \ref{39bus}. The assignments on parameters and settings are basically referred to \cite{chenxin,yifu}. Particularly, there are 10 generator buses located at buses 30-39, and the others are load buses approximately serving 6 GW, where buses 12-20 are controllable. The objective function of load adjustment is $$f_i(d_i)=a_id_i^2+b_i\|d_i+c_i\|_1+f_{i,1}(d_i),$$ where for load buses 12-17, $a_i$ and $b_i$ are randomly selected in [0.5,2.5] p.u. and [1,1.5] p.u., respectively. For other controllable loads, $a_i$ and $b_i$ are randomly selected in [2.5,3] p.u. and [1.5,2] p.u., respectively. $c_i=0.15*i$. The amplitude of controllable load is $|d_i|\leq$ 1.5 p.u.. We assign a uniform inertia $M_i=$ 8 p.u. and damping coefficient $D_i=$ 1 p.u. for generator buses and all buses, respectively.
\begin{figure}[t]
   \centering
   \includegraphics[scale=0.5]{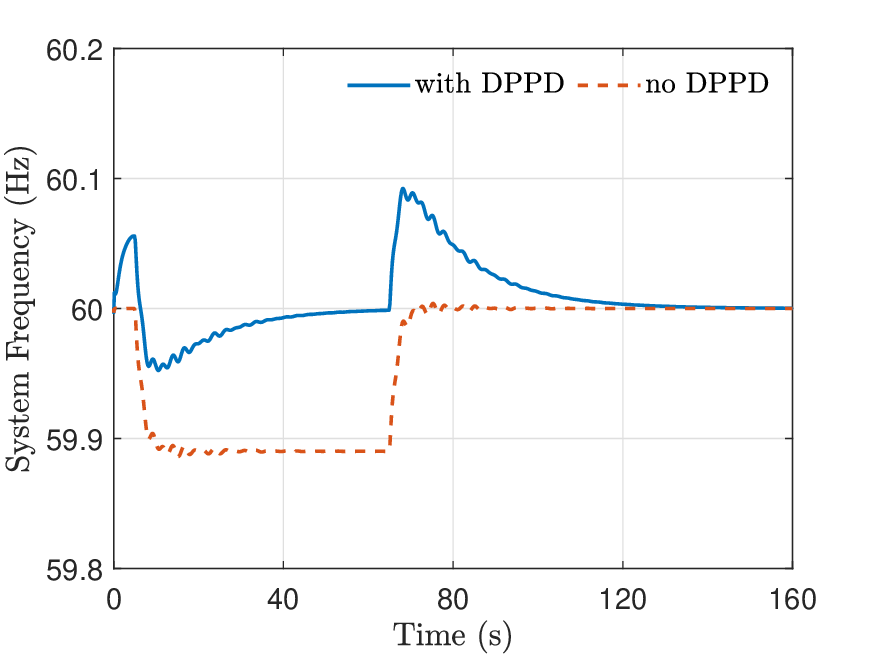}
   \caption{Frequency trajectories under step power injection change}\label{case1sysfre}
 \end{figure}

\subsection{Step power injection change}
In the first scenario, we consider constant power injections with the process of crash and recovery at certain buses. The comparison of frequency dynamics with and without the DPPD algorithm is illustrated in Fig. \ref{case1sysfre}. Initially, the power injections are set to be a balanced state such that the initial system frequency is 60 Hz. At time $t=$ 5s, a crash occurs at buses 37 and 39, where the power injections of generator buses 37 and 39 are specifically set to zero. Until $t=$ 65s, these two generators reconnect to the power network in the same states as before $t=$ 5s. It is obvious that
  \begin{figure}[t]
   \centering
   \includegraphics[scale=0.5]{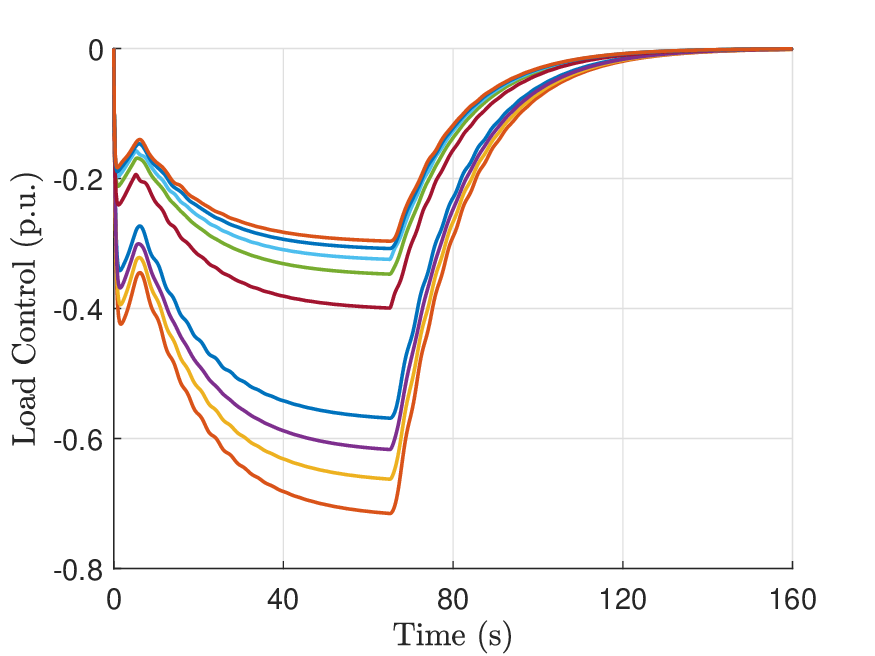}
   \caption{Load control inputs at buses 12-20 with DPPD algorithm}\label{case1loadcontrol}
 \end{figure}
the power system dynamics without DPPD algorithm are unable to restore the system frequency to its nominal value under imbalanced power injections. In contrast, the proposed DPPD algorithm facilitates the recovery process. Fig. \ref{case1loadcontrol} illustrates the evolution of control input at each controllable load bus. It shows that the load control inputs converge to zero from time $t=$ 65s for the rebalanced power injection. Furthermore, one can see that a lower cost coefficient corresponds to a larger control input.

\subsection{Time-varying power injection}
\begin{figure}[t]
    \subfigure[Initially balanced state]{
    \begin{minipage}[t]{0.48\linewidth}
    \includegraphics[scale=0.31]{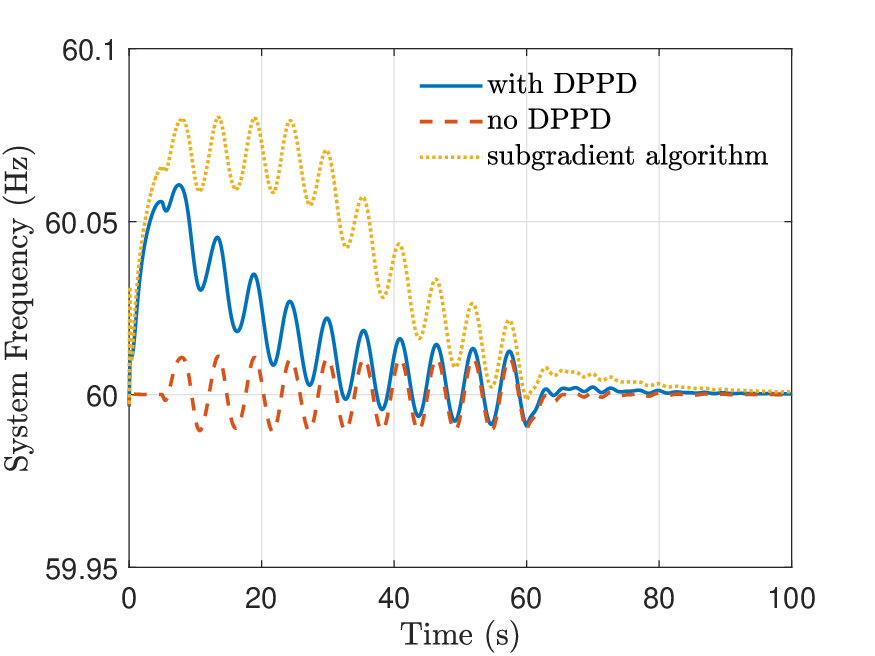}
    \label{case2sysfre_ba}
    \end{minipage}%
    }%
    \subfigure[Initially imbalanced state]{
    \begin{minipage}[t]{0.48\linewidth}
    \includegraphics[scale=0.31]{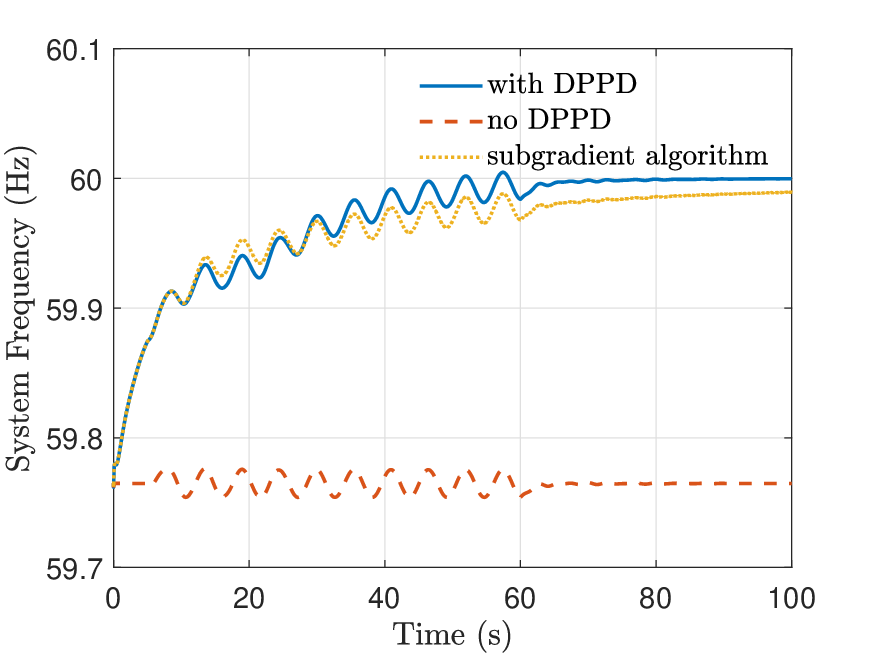}
    \label{case2sysfre_un}
    \end{minipage}
    }\caption{Frequency trajectories under time-varying power injections}\label{case2sysfre}
\end{figure}

\begin{figure}[t]
    \subfigure[Initially balanced state]{
    \begin{minipage}[t]{0.48\linewidth}
    \includegraphics[scale=0.31]{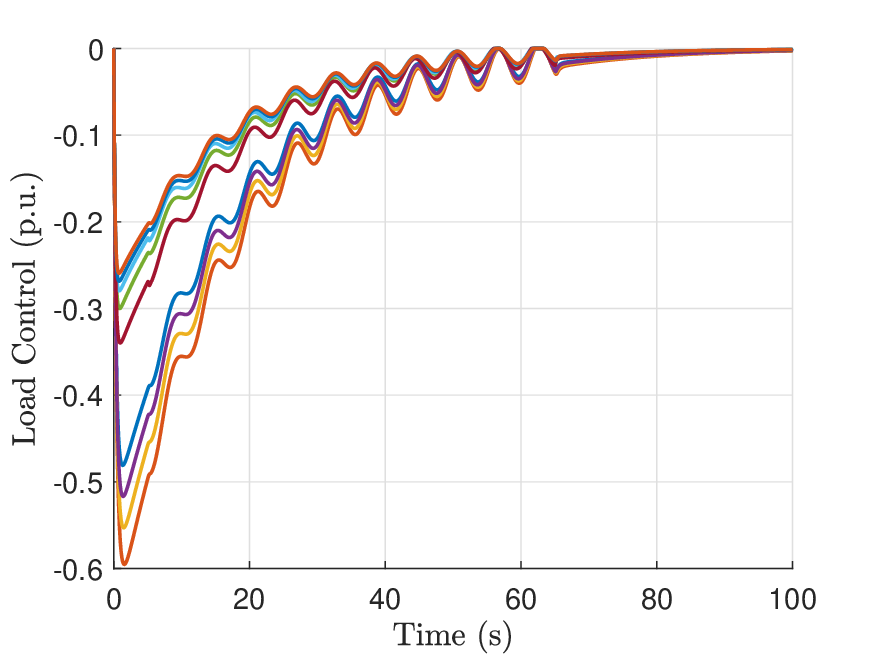}
    \label{case2loadcontrol_ba}
    \end{minipage}%
    }%
    \subfigure[Initially imbalanced state]{
    \begin{minipage}[t]{0.48\linewidth}
    \includegraphics[scale=0.31]{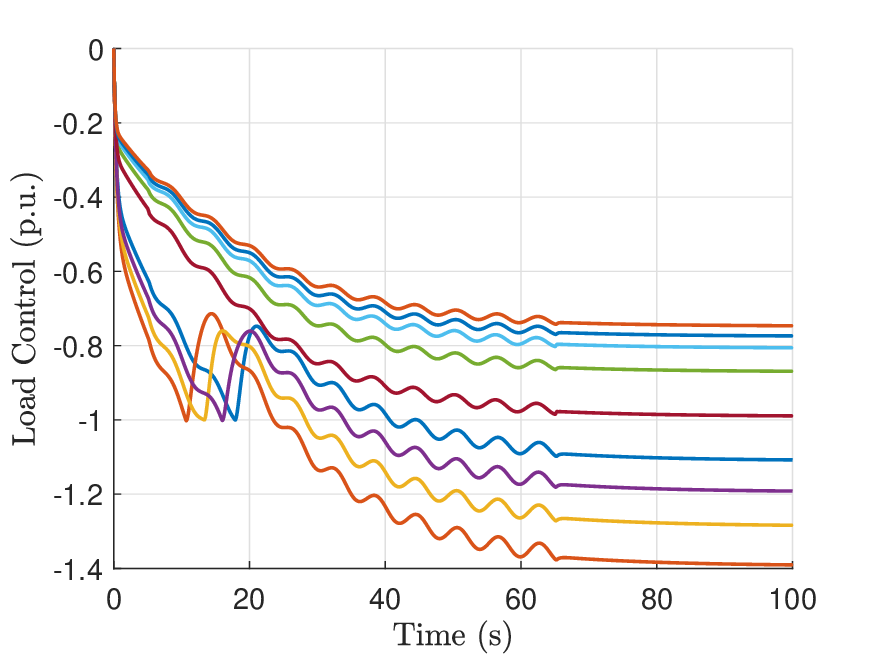}
    \label{case2loadcontrol_un}
    \end{minipage}
    }\caption{Load control inputs at buses 12-20 with DPPD algorithm}\label{case2loadcontrol}
\end{figure}

The performance of DPPD scheme under time-varying power injections is further studied, Compared with the power system dynamics without DPPD and automatic load control (ALC) scheme in \cite{chenxin} that employs with subgradient algorithm. In this scenario, the power injections at buses 37 and 39 are perturbed with a sinusoidal function that is proportional to their initial values during the time interval $[5,65]$s, and given as
\begin{equation*}
 P_i^{in}(t) = (1+0.4\sin(\frac{\pi t}{3}))P_i^{in}(0).
\end{equation*}
The remaining buses maintain their initial power injections throughout. Under such time-varying power injections, Fig. \ref{case2sysfre_ba} and Fig. \ref{case2sysfre_un} display the system frequency response with initially balanced and imbalanced power injections, respectively. It is seen that DPPD algorithm can maintain the frequency around 60 Hz in both initially balanced and imbalanced power injection scenarios, owning to the feedback of frequency deviation and local power imbalance in the designed DPPD algorithm. Conversely, the distributed frequency controller without DPPD algorithm is effective only in initially balanced case. Regarding ALC scheme equipped with subgradient algorithm (referred to as subgradient algorithm for short in Fig. \ref{case2sysfre}), the trajectory of system frequency exhibits similar transient performance as applying DPPD, but the convergence rate is relatively slow. Fig. \ref{case2loadcontrol_ba} and Fig. \ref{case2loadcontrol_un} present the evolution of load control inputs with DPPD algorithm under perturbed power injections. As a comparison, the evolution of load control inputs using subgradient-based ALC scheme is tested in Fig. \ref{case2loadcontrol} with initially balanced power injection. From Fig. \ref{case2loadcontrol}, it reveals that the trajectories of some buses go zigzag at their nonsmooth points, i.e., $d_i=-0.15i$ for $i\in\{1,2,3\}$, owning to the direct employment of subgradient in the controller. In certain cases, persistent nonsmooth control is actually deemed unacceptable for controllers in physical systems.

\begin{figure}[htbp]
   \centering
   \includegraphics[scale=0.5]{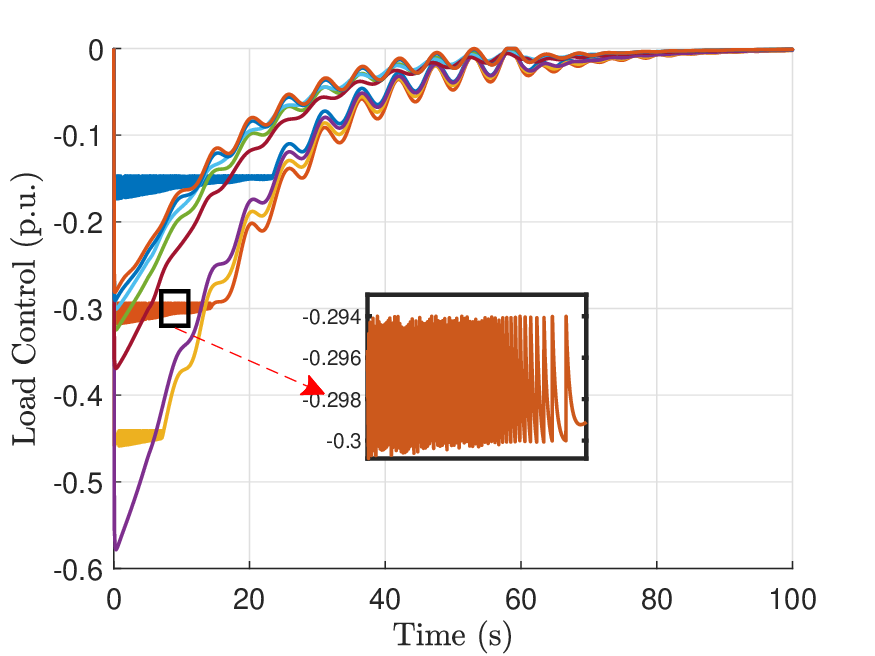}
   \caption{Load control inputs at buses 12-20 with subgradient-based ALC scheme}\label{case2loadcontrol}
\end{figure}

\subsection{Robustness against uncertainties}
\begin{figure}[t]
    \subfigure[$k_1=0.15$]{
    \begin{minipage}[t]{0.48\linewidth}
    \includegraphics[scale=0.31]{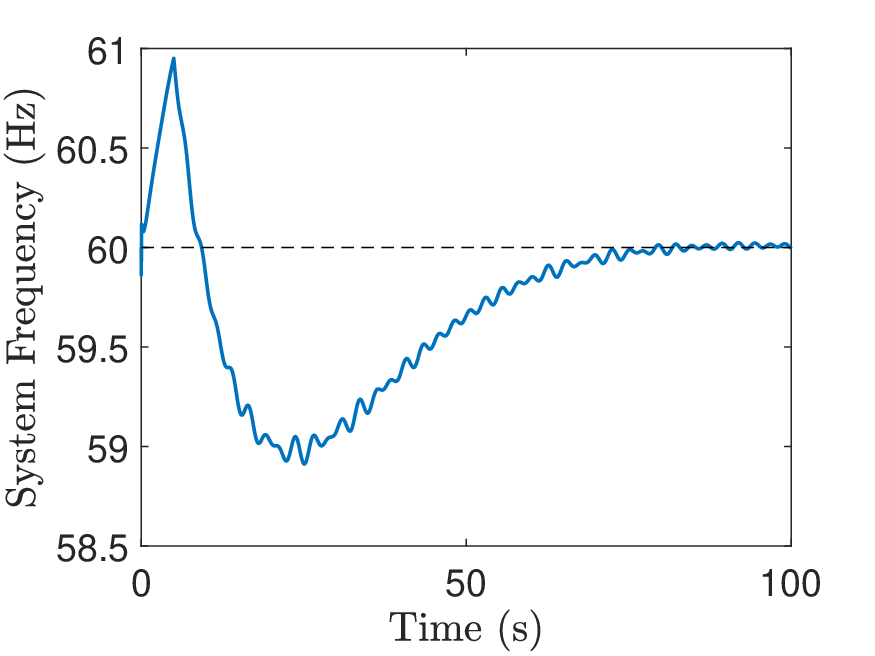}
    \label{case3a}
    \end{minipage}%
    }%
    \subfigure[$k_1=0.5$]{
    \begin{minipage}[t]{0.48\linewidth}
    \includegraphics[scale=0.31]{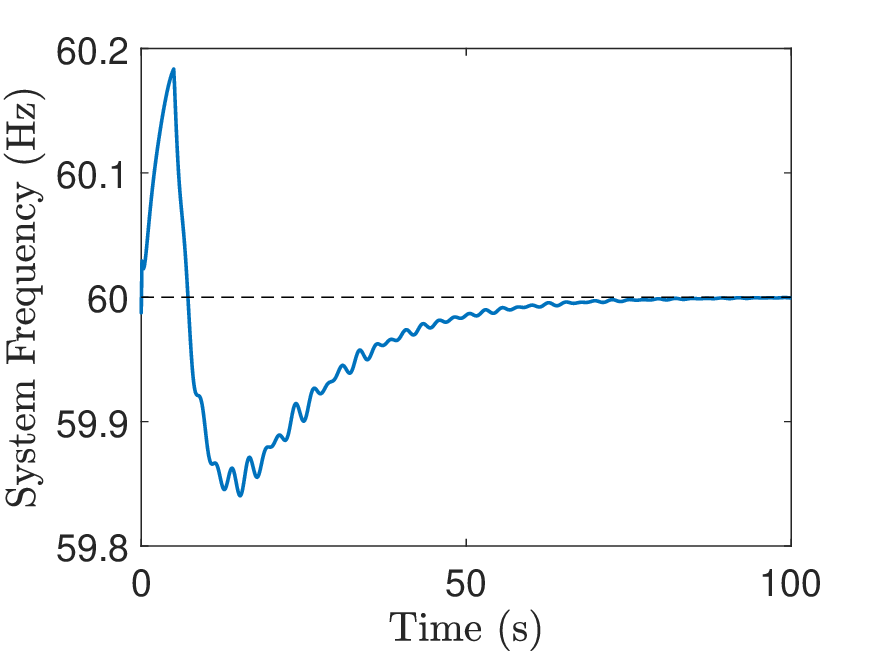}
    \label{case3b}
    \end{minipage}
    }%
\quad
    \subfigure[$k_1=2$]{
    \begin{minipage}[t]{0.48\linewidth}
    \includegraphics[scale=0.31]{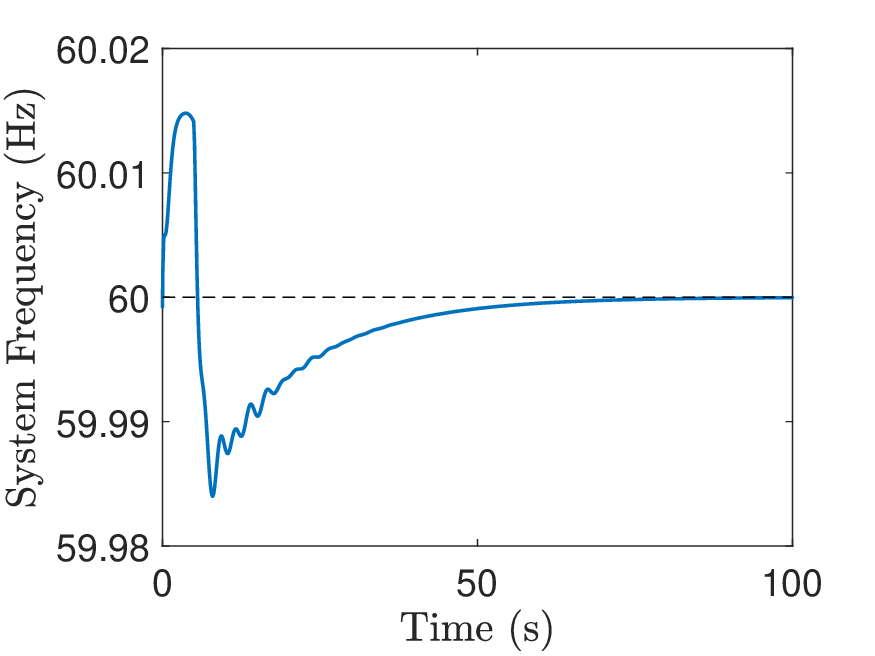}
    \label{case3c}
    \end{minipage}%
    }%
    \subfigure[$k_1=10$]{
    \begin{minipage}[t]{0.48\linewidth}
    \includegraphics[scale=0.31]{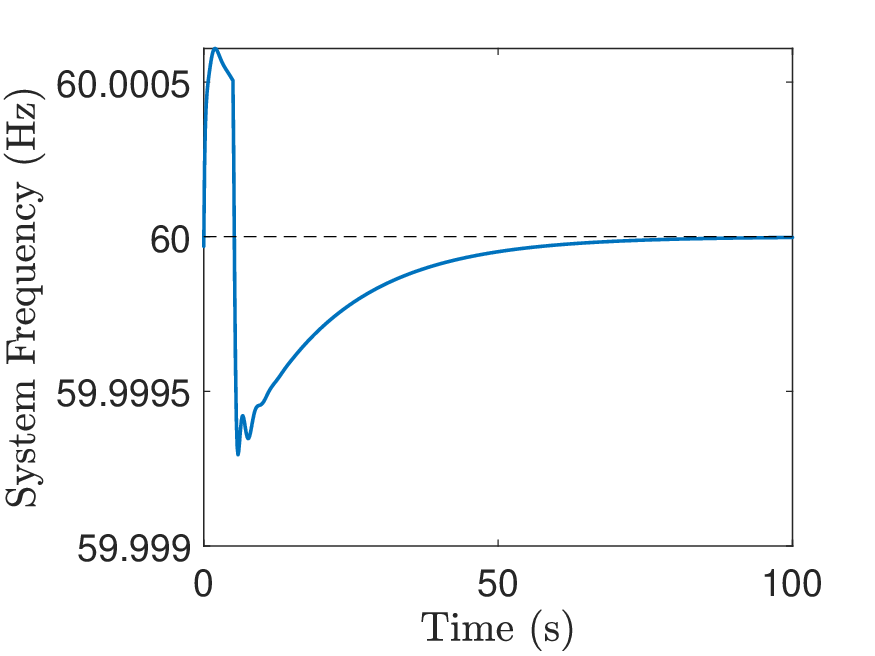}
    \label{case3d}
    \end{minipage}
    }%
    \caption{Frequency trajectories with different damping coefficients}\label{case3abcd}
\end{figure}

In this case study, it aims to illustrate the robustness of DPPD against two types of uncertainties, both of which are tested under a step power injection change. The first type of uncertainty pertains to non-ideal damping coefficient, where the damping coefficient is $k_1$ times the accurate one. The responses of system frequency with different damping coefficients for all buses are shown in Fig. \ref{case3abcd}. It indicates that a larger damping coefficient results in a tighter oscillation region of frequency, and consequently, leads to faster convergence towards the nominal frequency. This is because it requires a smaller response of frequency to counteract the power imbalance under a larger damping coefficient.

\begin{figure}[t]
    \subfigure[$k_2=0.15$]{
    \begin{minipage}[t]{0.48\linewidth}
    \includegraphics[scale=0.31]{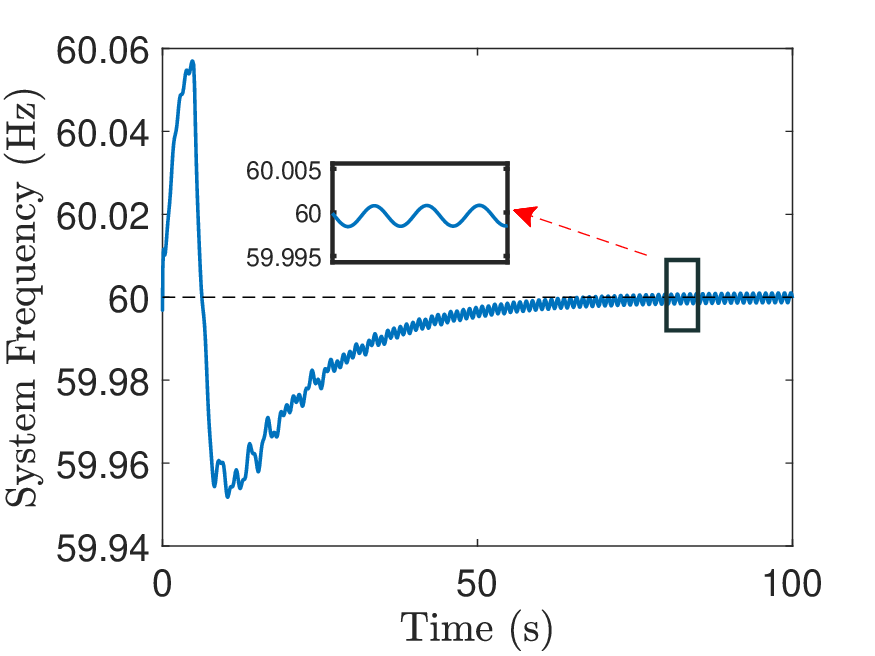}
    \label{case32a}
    \end{minipage}%
    }%
    \subfigure[$k_2=0.5$]{
    \begin{minipage}[t]{0.48\linewidth}
    \includegraphics[scale=0.31]{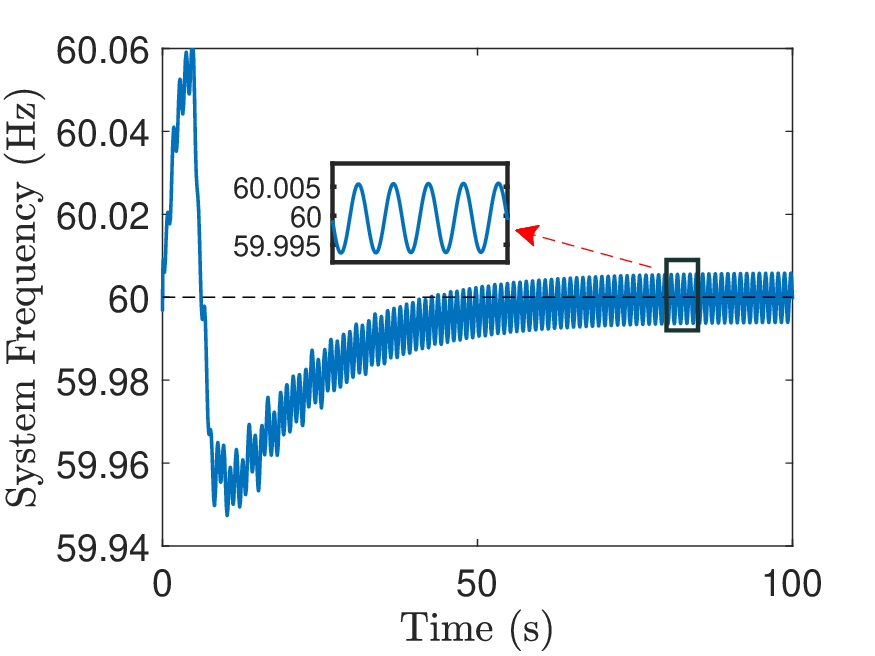}
    \label{case32b}
    \end{minipage}
    }\caption{Frequency trajectories with measurement errors in frequency}\label{case3omegnoise}
\end{figure}

The second type of uncertainty arises from measurement errors in $\omega$. In this case, the frequency measured at the controllable load bus is represented as
$$\omega_i^{'}(t)=\omega_i(t)+k_2\sin(2\pi t),\quad  i\in\{12,...,20\},$$
where $k_2$ is a tuned parameter. As shown in Fig. \ref{case3omegnoise}, it reflects that a greater amplitude of measurement noise in the local frequency results in a larger oscillation of the steady-state frequency. Generally, our proposed DPPD algorithm is capable to bring the system frequency back to the nominal value under tolerable measurement noise in frequency.

\section{Conclusion}
This paper has investigated an optimal frequency regulation issue on load side in power systems with nonsmooth cost functions. A distributed algorithm named DPPD was designed to solve such nonsmooth OLC problem, which can be deemed as an extension as PD-type (or PI-type) protocol compared to prior works. The optimality and global asymptotic convergence of DPPD algorithm were proved with a convergence rate as $O(1/\sqrt{t})$. Finally, simulations on the IEEE 39-bus system confirmed that the proposed DPPD algorithm is effective to restore the system frequency or maintain the system frequency within relatively safe region in the transients. Additionally, it is also robust to uncertainties of damping coefficients and measurement noises in frequency. Currently, limited by the theoretical analysis framework, we only establish a sublinear rate. Our future work will exploit linear convergence rate and improve the performance of transient frequency.

\section*{Appendix}

\subsection{Proof of Lemma 2}
\begin{IEEEproof}
    Denote an optimal solution of ROLC problem (3) as $(d^*,\hat{\theta}^*,\omega^*,P^*)$. First, according to the optimality condition of an optimal solution to ROLC problem in Lemma 1, and the obtained conclusions in the proof of Theorem 1, it gives $\omega_i^*=0, \forall i\in\mathcal{N}$.
    Then, let $\hat{\theta}^*=\theta^*+\epsilon\mathbf{1}$ and $P^*=BC^T\theta^*$ where $\epsilon$ is a constant. By comparing the formulations of problem (2) and problem (3), it follows that $(d^*,\theta^*,\omega^*=\mathbf{0})$ is a feasible solution of problem (2).
    Suppose that $(d^*,\theta^*,\omega^*=\mathbf{0})$ is not optimal with regard to problem (2), then it means there exists another feasible solution $(\tilde{d}^*,\tilde{\theta}^*,\tilde{\omega}^*)$ satisfying $\sum_{i\in\mathcal{N}}f_i(\tilde{d}_i^*)<\sum_{i\in\mathcal{N}}f_i(d_i^*)$. According to the formulation of problem (2), it must have $\tilde{\omega}^*=\mathbf{0}$, leading to $$\sum_ {i\in\mathcal{N}}f_i(\tilde{d}_i^*)+\frac{D_i}{2}\tilde{\omega}_i^{*2}<\sum_{i\in\mathcal{N}}f_i(d_i^*)+\frac{D_i}{2}\omega_i^{*2}.$$
    That is to say, by letting $\tilde{\theta}^*=\hat{\theta}^*$, it follows that $(\tilde{d}^*,\tilde{\theta}^*,\tilde{\omega}^*,P^*)$ is a feasible solution to problem (3) with a strictly lower cost than $(d^*,\hat{\theta}^*,\omega^*,P^*)$. However, it contradicts with the fact that $(d^*,\hat{\theta}^*,\omega^*,P^*)$ is optimal to problem (3). Thus, $(d^*,\theta^*,\omega^*)$ is exactly an optimal solution of problem (2) if $C^T\theta^*=C^T\hat{\theta}^*$ (i.e., $\hat{\theta}^*=\theta^*+\epsilon\mathbf{1}$) and $BC^T\theta^*=P^*$. This completes the proof.
\end{IEEEproof}

\subsection{Proof of Theorem 2}
\begin{IEEEproof}
 1) To show statement 1, let us see the distance of projecting any point $d$ onto the feasible region. Define a Lyapunov function candidate as
\begin{equation*}
  V_a= \frac{1}{2}\|\mathbf{prox}_{f_1}(d)-d\|^2,
\end{equation*}
where $\mathbf{prox}_{f_1}(d):=\mathrm{col}\{\mathbf{prox}_{f_{1,1}}(d_1),...,\mathbf{prox}_{f_{n,1}}(d_n)\}$. The derivative of $V_a$ can be written as
\begin{eqnarray*}
   \hspace{-17pt}\dot{V}_a\hspace{-7pt}&=&\hspace{-7pt}-(\mathbf{prox}_{f_1}(d)-d)^T\dot{d}\\
  \hspace{-17pt} \hspace{-7pt}&=&\hspace{-7pt}(\mathbf{prox}_{f_1}(d)-d)^T(\mathbf{prox}_{f_1}(d)-\mathbf{prox}_{f_1}^{d})\\
  \hspace{-17pt} \hspace{-7pt}&&\hspace{-7pt}-\|\mathbf{prox}_{f_1}(d)-d\|^2\\
\hspace{-17pt}\hspace{-7pt}&\leq&\hspace{-7pt}-\|\mathbf{prox}_{f_1}(d)-d\|^2,
\end{eqnarray*}
where $\mathbf{prox}_{f_{1}}^{d}:=\mathbf{prox}_{f_{1}}(d-\nabla f_0(d)+\kappa \eta +\lambda+u_1+\mu+u_2)$. The second equality is obtained by adding and subtracting $\mathbf{prox}_{f_1}(d)$ in the dynamics of $\dot{d}$. The inequality holds for the projection property on convex analysis, i.e., $$(\mathbf{prox}_{f_1}(d)-d)^T(\mathbf{prox}_{f_1}(d)-\mathbf{prox}_{f_1}^{d})\leq 0.$$

In light of $\dot{V}_a\leq 0$, we can claim that $d$ will always  remain within the power limit constraint if $d(0)\in \Omega$.
This completes the proof of the first statement.

 2) What follows is the proof of statement 2. For the analysis of this part, we equivalently rewrite \eqref{dppdd}, \eqref{dppdf} and \eqref{dppdg} as
 \begin{equation}\label{fuzhuzhengming}
    \begin{aligned}
   \dot{\hat{\theta}}_i&\hspace{-1.5pt}=\hspace{-1.5pt}\mathcal{P}_{\Theta_i}\hspace{-1.5pt}\Big [\hat{\theta}_i
   \hspace{-1.5pt}+\hspace{-5pt}\sum_{{j:ij\in\mathcal{E}}}\hspace{-5pt}B_{ij}(\mathcal{P}_{\mathbb{R}}^{\mu_i}\hspace{-1.5pt}-\hspace{-1.5pt}\mathcal{P}_{\mathbb{R}}^{\mu_j})
    \hspace{-1.5pt}-\hspace{-3.5pt}\sum_{{k:ki\in\mathcal{E}}}\hspace{-6.5pt}B_{ki}(\mathcal{P}_{\mathbb{R}}^{\mu_k}\hspace{-1.5pt}-\hspace{-1.5pt}\mathcal{P}_{\mathbb{R}}^{\mu_i})\Big]\hspace{-2pt}-\hspace{-2pt}\hat{\theta}_i,\\
   \dot{\omega}_i&\hspace{-1.5pt}=\hspace{-1.5pt}\mathcal{P}_{\mathbb{R}}^{\omega_i}-\omega_i,\\
   \dot{\mu}_i&\hspace{-1.5pt}=\hspace{-1.5pt}\mathcal{P}_{\mathbb{R}}^{\mu_i}-\mu_i,
    \end{aligned}
 \end{equation}
  where $\mathcal{P}_{\mathbb{R}}^{\omega_i}:=\mathcal{P}_{\mathbb{R}}(\omega_i+u_{i,1})$ and $\mathcal{P}_{\mathbb{R}}^{\mu_i}:=\mathcal{P}_{\mathbb{R}}(\mu_i+u_{i,2})$. $\Theta_i$ denotes the local feasible set of virtual phase angle, and $\Theta=\Theta_1\times\ldots\times\Theta_n$ being their Cartesian product. Note that it is a more general formulation such that our analysis is also applicable to the case where the operational constraint $\Theta$ is a general convex set.
 Besides, denote $\Delta d=d-d^*$, and doing likewise with $\eta$, $\hat{\theta}$, $P$, $\omega$ and $\mu$.

  Let construct a Lyapunov function candidate as $$V_b=V_1+V_2+V_3+V_4,$$ whose details will be unfolded in the sequel.

 Firstly, according to the gradient dynamics in the DPPD algorithm, define an auxiliary function as follows: $$\Psi(d,\hat{\theta},P,\omega,\mu)=f_0(d)+\frac{1}{2}\|\mathcal{P}_{\mathbb{R}^n}^{\omega}\|^2+\frac{1}{2}\|\mathcal{P}_{\mathbb{R}^n}^{\mu}\|^2.$$ Meanwhile, let
\begin{eqnarray*}
  V_1 \hspace{-6pt}&=& \hspace{-6pt} \Psi-\Psi^*-\nabla_{d}^T{\Psi^*}\Delta d-{\omega^*}^T(I-D)\Delta\omega-{\mu^*}^T\Delta\mu\\
  &&\hspace{-6pt}+{\omega^*}^TC\Delta P+{\mu^*}^TCBC^T\Delta \hat{\theta},
\end{eqnarray*}
where $\nabla_{d}{\Psi^*}$ denotes $\nabla_{d}{\Psi}({d}^*,{\hat{\theta}}^*,P^*,\omega^*,\mu^*)$, $\mathcal{P}_{\mathbb{R}^n}^{\omega}:=\mathrm{col}\{\mathcal{P}_{\mathbb{R}}^{\omega_1},...,\mathcal{P}_{\mathbb{R}}^{\omega_n}\}$ and $\mathcal{P}_{\mathbb{R}^n}^{\mu}:=\mathrm{col}\{\mathcal{P}_{\mathbb{R}}^{\mu_1},...,\mathcal{P}_{\mathbb{R}}^{\mu_n}\}$. Note that $f_0(d)$, $\mathcal{P}_{\mathbb{R}^n}^{\omega}$, $\mathcal{P}_{\mathbb{R}^n}^{\mu}$ and the norm operator $\|\cdot\|^2$ are all convex. As a consequence, it is straightforward to see the convexity of $\Psi$, together with $V_1\geq 0$ due to the property of convex function. Therewith, one can calculate the derivative of $V_1$ as
\begin{eqnarray*}
  \dot{V}_1\hspace{-6pt}&=&\hspace{-6pt}\underbrace{(\nabla_{d} \Psi-\nabla_{d}{\Psi^*})^T\dot{d}}_{S_{1.1}}+\underbrace{(\nabla_\omega\Psi-(I-D)\omega^*)^T\dot{\omega}}_{S_{1.2}}\nonumber\\
   \hspace{-6pt}&&\hspace{-6pt}+\underbrace{(\nabla_\mu\Psi-\mu^*)^T\dot{\mu}}_{S_{1.3}}+\underbrace{(\nabla_{\hat{\theta}}\Psi+CBC^T{\mu^*})^T\dot{\hat{\theta}}}_{S_{1.4}}\nonumber\\
   \hspace{-6pt}&&\hspace{-6pt}+\underbrace{(\nabla_{P}\Psi+C^T\omega^*)^T\dot{P}}_{S_{1.5}},
\end{eqnarray*}
where, for the subsequent analysis, $S_{1.2}$ is further rewritten as $S_{1.2}=S_{1.21}+S_{1.22}$ with
\begin{equation*}
\begin{aligned}
   S_{1.21}&=(\mathcal{P}_{\mathbb{R}^n}^\omega-\omega^*)^T\dot{\omega},  \\
   S_{1.22}&=-(\mathcal{P}_{\mathbb{R}^n}^\omega-\omega^*)^TD\dot{\omega}.
\end{aligned}
\end{equation*}

Now, define
$$V_2=\frac{1}{2}(\|\Delta {d}\|^2-2\kappa\Delta {d}^T\Delta\eta+\kappa\|\Delta\eta\|^2).$$
It is direct to verify the positive definiteness of $V_2$ under the condition $\kappa\in(0,1)$. With $\zeta:=\dot{V}_2+S_{1.1}$, one has
\begin{equation}\label{zeta}
\begin{aligned}
   &\zeta= \Delta {d}^T\dot{d}+\kappa\Delta\eta^T\dot{\eta}-\kappa\Delta\eta^T\dot{d}-\kappa\Delta {d}^T\dot{\eta}\nonumber\\
   &\quad\,\,\,\,\,+(\nabla f_0(d)-\nabla f_0({d^{*}}))^T\dot{d}-(\omega+u_1)^T\dot{d}\nonumber\\
   &\quad\,\,\,\,\,+(\omega^*+u_1^*)^T\dot{d}-(\mu+u_2)^T\dot{d}+(\mu^*+u_2^*)^T\dot{d}.
\end{aligned}
\end{equation}
Meanwhile, from dynamics \eqref{dppda} and \eqref{dppdb}, we can sort out the following equations
\begin{eqnarray}\label{sige}
  -\kappa\eta\hspace{-8pt}&=&\hspace{-8pt}\dot{\eta}+\partial f_2(\dot{\eta}+d),\nonumber\\
  -\kappa\eta^*\hspace{-8pt}&=&\hspace{-8pt}\partial f_2({d^{*}}),\nonumber\\
  -\nabla f_0({d})+\kappa\eta+\omega+\mu+u\hspace{-8pt}&=&\hspace{-8pt}\dot{d}+\partial f_1(d+\dot{d}), \nonumber\\
  -\nabla f_0({d^{*}})+\kappa\eta^*+\omega^*+\mu^*+u^*\hspace{-8pt}&=&\hspace{-8pt}\partial f_1({d^{*}}),
\end{eqnarray}
where, for brevity, $u=u_1+u_2$ and $u^*=u_1^*+u_2^*$.
On account of the maximal monotonicity of $\partial f_1(\cdot)$ and $\partial f_2(\cdot)$, it follows from \eqref{sige} that
\begin{eqnarray}
   &&\hspace{-40pt} \langle-\kappa\Delta\eta-\dot{\eta},\Delta d+\dot{\eta}\rangle\geq 0,\nonumber\\
   &&\hspace{-40pt}\langle-\nabla f_0(d)+\nabla f_0({d^{*}})+\Delta \omega+\Delta\mu\nonumber\\
    &&\hspace{-40pt}\qquad\qquad+u-u^*+\kappa\Delta\eta-\dot{d},\Delta d+\dot{d}\rangle\geq0.
\end{eqnarray}
By rearranging the above two inequalities, and substituting \eqref{zeta} into here, it boils down to
\begin{eqnarray}\label{16}
   \zeta\hspace{-7pt}&\leq&\hspace{-7pt} -\|\dot{d}\|^2-\|\dot{\eta}\|^2-(\nabla f_0(d)-\nabla f_0({d^{*}}))^T\Delta d\nonumber\\
  \hspace{-7pt}&&\hspace{-7pt} -(1+\kappa)\Delta {d}^T\dot{\eta}+(\Delta \omega+u_1-u_1^*)^T\Delta d \nonumber\\
   \hspace{-7pt}&&\hspace{-7pt}+(\Delta \mu+u_2-u_2^*)^T\Delta d.
\end{eqnarray}
Under Assumption 1, it further has
\begin{eqnarray}\label{17}
    \hspace{-25pt}&&-(1+\kappa)\Delta {d}^T\dot{\eta}\nonumber\\
    \hspace{-25pt}&&\leq\frac{1+\kappa}{2\alpha}\|\Delta d\|^2 +\frac{\alpha(1+\kappa)}{2}\|\dot{\eta}\|^2\nonumber \\
    \hspace{-25pt}&&\leq\frac{1+\kappa}{2\alpha\beta}(\nabla f_0(d)-\nabla f_0({d^{*}})^T\Delta d+\frac{\alpha(1+\kappa)}{2}\|\dot{\eta}\|^2,
\end{eqnarray}
where $\alpha$ is a positive constant.
Then, we can deduce from \eqref{16} and \eqref{17} that
\begin{eqnarray*}
   \zeta\hspace{-8pt}&\leq&\hspace{-8pt} -\|\dot{d}\|^2-(1-\frac{1+\kappa}{2\alpha\beta})(\nabla f_0(d)-\nabla f_0({d^{*}}))^T\Delta d \nonumber\\
   &&\hspace{-8pt} -(1-\frac{\alpha(1+\kappa)}{2})\|\dot{\eta}\|^2+\underbrace{(\Delta \omega+u_1-u_1^*)^T\Delta d}_{S_{2.1}}\nonumber\\
   &&\hspace{-8pt} +\underbrace{(\Delta \mu+u_2-u_2^*)^T\Delta d}_{S_{2.2}}.
\end{eqnarray*}
At this stage, to illustrate the existence of $\alpha$ such that $1-\frac{1+\kappa}{2\alpha\beta}\geq 0$ and $1-\frac{\alpha(1+\kappa)}{2}\geq 0$, i.e., $\frac{1+\kappa}{2\beta}\leq\alpha\leq\frac{2}{1+\kappa}$, we define $h(\kappa)=\frac{2}{1+\kappa}-\frac{1+\kappa}{2\beta}$. Since $\frac{\partial h}{\partial \kappa}=-\frac{2}{(1+\kappa)^2}-\frac{1}{2\beta}<0$ for all $\kappa\in(0,1)$, that means $h_{\min}=h(1)=1-\frac{1}{\beta}$. Clearly, $h(\kappa)\geq 0$ always holds under Assumption 1. As such, $1-\frac{1+\kappa}{2\alpha\beta}$ and $1-\frac{\alpha(1+\kappa)}{2}$ are guaranteed to be nonnegative.

Next, in order to eliminate $S_{1.21}$ and $S_{2.1}$ as well as $S_{1.3}$ and $S_{2.2}$, define $$V_3=V_{3.1}+V_{3.2}+V_{3.3},$$ where $V_{3.1}=\frac{1}{2}\|\Delta \mu\|^2$, $V_{3.2}=\frac{1}{2}\|\Delta \omega\|^2$ and $V_{3.3}=\frac{3}{2}\Delta\omega^TD\Delta\omega$.

By giving insight into the derivative of $V_{3.1}$ with $S_{1.3}$, one can obtain that
\begin{eqnarray}\label{mu_1st}
&&\hspace{-25pt}\dot{V}_{3.1}+S_{1.3}\nonumber\\
&&\hspace{-25pt}=\Delta \mu^T\dot{\mu}+(\mathcal{P}_{\mathbb{R}^n}^\mu-\mu^*)^T\dot{\mu}\nonumber\\
&&\hspace{-25pt}=\frac{1}{2}(\mathcal{P}_{\mathbb{R}^n}^\mu-\mu)^T(\Delta\mu+\mathcal{P}_{\mathbb{R}^n}^\mu-\mu^*) \nonumber\\
&&\hspace{-25pt}=-\frac{1}{2}\|\mathcal{P}_{\mathbb{R}^n}^\mu-\mu\|^2+(\mathcal{P}_{\mathbb{R}^n}^\mu-\mu)^T(\mathcal{P}_{\mathbb{R}^n}^\mu-\mu^*)\nonumber\\
&&\hspace{-25pt}=-\frac{1}{2}\|\mathcal{P}_{\mathbb{R}^n}^\mu-\mu\|^2-(\mathcal{P}_{\mathbb{R}^n}^\mu-\mu^*)^T(d+CBC^T\hat{\theta}-P^{in})\nonumber\\
&&\hspace{-25pt}\quad+(\mathcal{P}_{\mathbb{R}^n}^\mu-\mu^*)^T(\mathcal{P}_{\mathbb{R}^n}^\mu-\mu-P^{in}+d+CBC^T\hat{\theta})\nonumber\\
&&\hspace{-25pt}=-\frac{1}{2}\|\mathcal{P}_{\mathbb{R}^n}^\mu-\mu\|^2-(\mathcal{P}_{\mathbb{R}^n}^\mu-\mu^*)^T\Delta d\nonumber\\
&&\hspace{-25pt}\quad-(\mathcal{P}_{\mathbb{R}^n}^\mu-\mu^*)^TCBC^T\Delta\hat{\theta}\nonumber\\
&&\hspace{-25pt}\quad-(\mathcal{P}_{\mathbb{R}^n}^\mu-\mu^*)^T(d^{*}+CBC^T\hat{\theta}^*-P^{in})\nonumber\\
&&\hspace{-25pt}=-\frac{1}{2}\|\mathcal{P}_{\mathbb{R}^n}^\mu-\mu\|^2\underbrace{-(\mathcal{P}_{\mathbb{R}^n}^\mu-\mu^*)^T\Delta d}_{S_{3.1}}\nonumber\\
&&\hspace{-25pt}\quad\underbrace{-(\mathcal{P}_{\mathbb{R}^n}^\mu-\mu^*)^TCBC^T\Delta\hat{\theta}}_{S_{3.2}},
\end{eqnarray}
where, the facts that $(\mathcal{P}_{\mathbb{R}^n}^\mu-\mu^*)^T(\mathcal{P}_{\mathbb{R}^n}^\mu-\mu-P^{in}+d+CBC^T\hat{\theta})=0$ and $-(\mathcal{P}_{\mathbb{R}^n}^\mu-\mu^*)^T(d^{*}+CBC^T\hat{\theta}^*-P^{in})=0$ are used in the fourth and fifth equality, respectively. Whereupon, it has
\begin{equation}\label{23}
S_{2.2}+S_{3.1}=0.
\end{equation}
Similarly, the derivative of $V_{3.2}$ with $S_{1.21}$ is calculated as
\begin{eqnarray}\label{omega_1st}
    &&\hspace{-22pt}\dot{V}_{3.2}+S_{1.21}\nonumber\\
    &&\hspace{-22pt}=-\frac{1}{2}\|\mathcal{P}_{\mathbb{R}^n}^\omega-\omega\|^2+(\mathcal{P}_{\mathbb{R}^n}^\omega-\omega^*)^T(\mathcal{P}_{\mathbb{R}^n}^\omega-\omega)\nonumber\\
    &&\hspace{-22pt}=-\frac{1}{2}\|\mathcal{P}_{\mathbb{R}^n}^\omega-\omega\|^2\underbrace{-(\mathcal{P}_{\mathbb{R}^n}^\omega-\omega^*)^T\Delta d}_{S_{3.3}}\nonumber\\
    &&\hspace{-22pt}\quad\underbrace{-(\mathcal{P}_{\mathbb{R}^n}^\omega-\omega^*)^TD\omega}_{S_{3.4}}\underbrace{-(\mathcal{P}_{\mathbb{R}^n}^\omega-\omega^*)^TC\Delta P}_{S_{3.5}},
\end{eqnarray}
which yields
$$S_{2.1}+S_{3.3}=0.$$
Furthermore, $S_{1.22}$ and $S_{3.4}$ can be rearranged as
\begin{eqnarray*}
   &&\hspace{-20pt}S_{1.22}+S_{3.4}\nonumber\\
   &&\hspace{-20pt}=\hspace{-1.5pt}-(\mathcal{P}_{\mathbb{R}^n}^\omega\hspace{-1.2pt}-\hspace{-1.2pt}\omega)^TD\dot{\omega}\hspace{-1.2pt}-\hspace{-1.2pt}\Delta \omega^TD\dot{\omega}\hspace{-1.2pt}-\hspace{-1.2pt}(\mathcal{P}_{\mathbb{R}^n}^\omega\hspace{-1.2pt}-\hspace{-1.2pt}\omega)^TD\omega \hspace{-1.2pt}-\hspace{-1.2pt}\Delta\omega^TD\omega.
\end{eqnarray*}
 By combining with dynamics $\dot{\omega}_i=\mathcal{P}_{\mathbb{R}}^{\omega_i}-\omega_i$ as in \eqref{fuzhuzhengming}, the above equation implies that
\begin{eqnarray*}
   &&\hspace{-12pt}S_{1.22}+S_{3.4}\nonumber\\
   &&\hspace{-12pt}=-2\dot{\omega}^TD\dot{\omega}-\Delta \omega^TD\dot{\omega}-2\dot{\omega}^TD\omega-\Delta\omega^TD\omega\nonumber\\
   &&\hspace{-12pt}\leq-\Delta\omega^TD\dot{\omega}-2\dot{\omega}^TD\omega\nonumber\\
   &&\hspace{-12pt}=-3\Delta \omega^TD\dot{\omega},
\end{eqnarray*}
where the inequality results from $-2\dot{\omega}^TD\dot{\omega}\leq 0$ and $-\Delta\omega^TD\omega\leq0$, due to the facts that $D$ is positive semidefinite and $\omega^*=\mathbf{0}$. With $V_{3.3}$, it leads to
\begin{equation*}
   \dot{V}_{3.3}+S_{1.22}+S_{3.4}\leq 0.
\end{equation*}

To proceed, in order to eliminate the residual terms, i.e., $S_{1.4}+S_{1.5}+S_{3.2}+S_{3.5}$, function $V_4$ is defined as
\begin{equation}\label{V4}
   V_{4}=\frac{1}{2}(\|\Delta\hat{\theta}\|^2+\|\Delta P\|^2).
\end{equation}
As the aforementioned proofs go, $V_4$ can also be divided into two groups, whose proofs are similar. The first group is about $\mu$ and $\hat{\theta}$, while the other group is about $\omega$ and $P$.

Firstly, for $\mu$ and $\hat{\theta}$, it can be calculated that
\begin{eqnarray}\label{theta_1st}
  &&\hspace{-8pt}\Delta \hat{\theta}^T\dot{\hat{\theta}}+S_{1.4}\nonumber \\
  &&\hspace{-8pt}=\Delta \hat{\theta}^T(\mathcal{P}_{\Theta}^{\hat{\theta}}-\hat{\theta})+(\nabla_{\hat{\theta}}\Psi-\nabla_{\hat{\theta}}\Psi^*)^T(\mathcal{P}_{\Theta}^{\hat{\theta}}-\hat{\theta})\nonumber\\
  &&\hspace{-8pt}=\Delta\hat{\theta}^T(\mathcal{P}_{\Theta}^{\hat{\theta}}-\hat{\theta})+(\nabla_{\hat{\theta}}\Psi-\nabla_{\hat{\theta}}\Psi^*)^T(\mathcal{P}_{\Theta}^{\hat{\theta}}-{\hat{\theta}}^*)\nonumber\\
  &&\hspace{-8pt}\quad-(\nabla_{\hat{\theta}}\Psi-\nabla_{\hat{\theta}}\Psi^*)^T\Delta\hat{\theta}\nonumber\\
  &&\hspace{-8pt}=\Delta\hat{\theta}^T(\mathcal{P}_{\Theta}^{\hat{\theta}}-\hat{\theta})+(\nabla_{\hat{\theta}}\Psi-\hat{\theta}+\mathcal{P}_{\Theta}^{\hat{\theta}})^T(\mathcal{P}_{\Theta}^{\hat{\theta}}-{\hat{\theta}}^*)\nonumber\\
  &&\hspace{-8pt}\quad-(\mathcal{P}_{\Theta}^{\hat{\theta}}-\hat{\theta})^T(\mathcal{P}_{\Theta}^{\hat{\theta}}-{\hat{\theta}}^*)-\nabla_{\hat{\theta}}^T\Psi^*(\mathcal{P}_{\Theta}^{\hat{\theta}}-{\hat{\theta}}^*)\nonumber\\
  &&\hspace{-8pt}\quad-(\nabla_{\hat{\theta}}\Psi-\nabla_{\hat{\theta}}\Psi^*)^T\Delta\hat{\theta}\nonumber\\
  &&\hspace{-8pt}\leq-\|\mathcal{P}_{\Theta}^{\hat{\theta}}-\hat{\theta}\|^2\underbrace{-(\nabla_{\hat{\theta}}\Psi-\nabla_{\hat{\theta}}\Psi^*)^T\Delta\hat{\theta}}_{S_{4.1
  }},
\end{eqnarray}
where $\mathcal{P}_{\Theta}^{\hat{\theta}}:=\mathcal{P}_{\Theta}(\hat{\theta}+CBC^T\mathcal{P}_{\mathbb{R}^n}^\mu)=\mathrm{col}\{\mathcal{P}_{\Theta_1}^{\hat{\theta}_1},...,\mathcal{P}_{\Theta_n}^{\hat{\theta}_n}\}$. To obtain the inequality, we have utilized the properties that $(\nabla_{\hat{\theta}}\Psi-\hat{\theta}+\mathcal{P}_{\Theta}^{\hat{\theta}})^T(\mathcal{P}_{\Theta}^{\hat{\theta}}-{\hat{\theta}}^*)\leq0$ and $-\nabla_{\hat{\theta}}^T\Psi^*(\mathcal{P}_{\Theta}^{\hat{\theta}}-{\hat{\theta}}^*)\leq 0$ in virtue of convex analysis.
As the similar way to obtaining \eqref{theta_1st} by convex analysis, it has the following result for those terms associated with $P$ and $\omega$, i.e.,
\begin{equation}
   \Delta P^T\dot{P}+S_{1.5}\leq-\|\mathcal{P}_{\mathbb{R}^l}^P-P\|^2\underbrace{-(\nabla_{P}\Psi-\nabla_{P}\Psi^*)^T\Delta P}_{S_{4.2
  }},
\end{equation}
where $\mathcal{P}_{\mathbb{R}^l}^{P}:=\mathcal{P}_{\mathbb{R}^l}[P+C^T\mathcal{P}_{\mathbb{R}^n}(\omega+u_1)]=\mathrm{col}\{\mathcal{P}_{\mathbb{R}}^{P_{ij}}\}_{ij\in\mathcal{E}}$. Then, with some calculations, it gives rise to

\begin{equation}\label{dixiao}
   S_{3.2}+S_{3.5}+S_{4.1}+S_{4.2}=0.
\end{equation}

Summarizing the above-discussed analysis, one can conclude that
\begin{eqnarray}\label{dot_Vb}
   \dot{V}_b\hspace{-8pt}&\leq&\hspace{-8pt} -\|\dot{d}\|^2-(1-\frac{1+\kappa}{2\alpha\beta})(\nabla f_0(d)-\nabla f_0(d^{*}))^T\Delta d \nonumber\\
   &&\hspace{-9pt} -(1-\hspace{-1pt}\frac{\alpha(1+\kappa)}{2})\|\dot{\eta}\|^2\hspace{-1pt}-\hspace{-1pt}\frac{1}{2}\|\mathcal{P}_{\mathbb{R}^n}^{\mu}-\mu\|^2\hspace{-1pt}-\hspace{-1pt}\frac{1}{2}\|\mathcal{P}_{\mathbb{R}^n}^{\omega}-\omega\|^2\nonumber\\
   &&\hspace{-9pt}-\|\mathcal{P}_{\Theta}^{\hat{\theta}}-\hat{\theta}\|^2-\|\mathcal{P}_{\mathbb{R}^l}^{P}-P\|^2,
\end{eqnarray}
which is negative semidefinite if Assumptions 1 and 2 hold, and $\kappa\in(0,1)$. At this moment, we have $V_b(t)\leq V_b(0)$. By observing the construction of $V_b$, it implies the boundedness of trajectory $(\eta,{d},{\hat{\theta}},P,\omega,\mu)$.


3) Now, in virtue of the LaSalle invariance principe \cite{lasalle}, one can conclude that the trajectory $(\eta,d,\hat{\theta},P,\omega,\mu)$ asymptotically converges to the largest invariant set contained in $\Xi$, where $\Xi=\{(\eta,d,\hat{\theta},P,\omega,\mu)\in\mathcal{M}:\dot{V}_b=0\}$ with $\mathcal{M}=\{(\eta,d,\hat{\theta},P,\omega,\mu)|V\leq q\}$, and $q>0$ is strongly positive invariant and bounded. To see the property of $\Xi$, one can observe from $\dot{V}_b$ that the point rendering $\dot{V}_b=0$ is exactly an equilibrium of \eqref{dppd}. Moreover, by Lemma 1 and Theorem 1, and noting that $V_b$ is radially unbounded, it follows that any trajectory is globally asymptotically convergent to an optimal solution of problem \eqref{rolc}.
This completes the proof.
\end{IEEEproof}


\end{document}